\magnification\magstep1
\baselineskip = 18pt
\def\n{\noindent}

\def\n{\noindent}
\def\vp{\varepsilon}
\def\cf{{\it cf.\/}\ }
\def\ent{{{\rm Z}\mkern-5.5mu{\rm Z}}}
 \overfullrule = 0pt
 \def\qed{{\hfill{\vrule height7pt width7pt
depth0pt}\par\bigskip}}

\let\eps =\varepsilon
\def \comp{ \;{}^{ {}_\vert }\!\!\!{\rm C}   }
\def\N#1{\left\Vert#1\right\Vert}
\def \nat{ { {\rm I}\!{\rm N}} }

\overfullrule = 0pt
\def\pf{\medskip{\noindent{\bf Proof. }}}
\centerline{\bf Bilinear forms on exact operator spaces}

\centerline{\bf and}

\centerline{\bf  $B(H)\otimes B(H)$.}
\centerline{by}
\centerline{Marius Junge and Gilles Pisier\footnote*{Partially
supported by
the NSF}}\vskip.5in

{\bf Abstract.}

Let $E,F$ be exact operator spaces (for example subspaces of the
$C^*$-algebra $K(H)$
of all the compact operators on an infinite dimensional Hilbert space
$H$). We
study a class of bounded linear maps $u\colon E\to F^*$ which we call
tracially
bounded. In particular, we prove that every completely bounded (in
short $c.b.$)
map $u\colon E\to F^*$ factors boundedly through a Hilbert space. This
is used to show that the set
$OS_n$ of all $n$-dimensional operator spaces equipped with the $c.b.$
version of
the Banach Mazur distance is not separable if $n>2$.

\n As an application we show that there is more than one $C^*$-norm on
$B(H)\otimes B(H)$, or equivalently that
$$B(H)\otimes_{\min}B(H)\not=B(H)\otimes_{\max}B(H),$$
which answers a long standing open question.

\n Finally we show that every ``maximal" operator space
(in the sense of Blecher-Paulsen)
is not exact in the infinite dimensional case, and in the finite
dimensional
case, we give a lower bound for the ``exactness constant".

\item{} Plan.

\item{0.} Introduction and background.

\item{1.} Factorization of bilinear forms on exact operator spaces.

\item{2.} The non-separability of $OS_n$.

\item{3.} Applications to $B(H)\otimes B(H)$ and maximal operator
spaces.

\item{4.} A new tensor product for $C^*$-algebras or operator spaces.

\vfill\eject

\n {\bf \S 0. Introduction and background.}

Following the remarkable results of Kirchberg on exact $C^*$-algebras
\cf
[Ki1, Ki2], the notion of ``exact operator space" was studied in the
paper [P1].
In this paper we continue the investigation started in [P1].

Let $E,F$ be operator spaces, we denote
$$d_{cb}(E,F)=\inf\{\N{u}_{cb}\N{u^{-1}}_{cb}\}$$
where the infimum runs over all possible isomorphisms $u\colon E\to F$.
If $E,F$
are not isomorphic, we set, by convention, $d_{cb}(E,F)=\infty$.
We denote by $B(H)$ (resp. $K(H)$)
the algebra of all bounded (resp. compact)
 operators on a Hilbert space $H$.
(See below for unexplained notation and terminology.)

In [P1], the following characteristic of an operator space, $E$ with
$\dim E=n$, was studied:
$$d_{SK}(E)=\inf\{d_{cb}(E,F)|\ F\subset K(\ell_2)\},\leqno(0.1)$$
or equivalently
$$d_{SK}(E)=\inf\{d_{cb}(E,F)|\ F\subset M_N,\ N\geq n\}.\leqno(0.1)'$$
For an infinite dimensional operator space $X$, we define
$$d_{SK}(X)=\sup\{d_{SK}(E)\}$$
where the supremum runs over all possible finite dimensional subspaces
$E\subset X$.

We say that $X$ is exact of $d_{SK}(X)<\infty$. This is equivalent (see
[P1])
to the exactness of a certain sequence of morphisms in the category of
operator
spaces, whence the terminology.

Let $E,F$ be exact operator spaces. We will prove below an inequality
for
completely bounded linear maps $u\colon E\to F^*$.
Here $F^*$ means the ``operator
space dual" of $F$ in the sense of [ER1, BP].

This key inequality can be viewed  as  a form of Grothendieck's
inequality for exact
operator spaces. It implies that every $c.b.$ map $u\colon E\to F^*$
can be factorized
(in the Banach space sense) through a Hilbert space with a norm of
factorization $\leq 4\|u\|_{cb}$. Actually, in this inequality
``complete
boundedness" can be replaced by a more general notion which we call
``tracial
boundedness" which has already been considered by previous authors (\cf
[I, B3]).
Using this key inequality and a somewhat surprising application of
Baire's
theorem, we prove that the metric space $OS_n$ of all $n$-dimensional
operator
spaces equipped with the distance
$$\delta_{cb}(E,F)=\log d_{cb}(E,F)$$
is not separable as soon as $n>2$. We can even show
that the subset of all isometrically
Hilbertian and homogeneous (in the sense of [P3])
$n$-dimensional operator spaces, is non-separable.

This has a surprising application to $C^*$-algebra
theory. Recall that a $C^*$-algebra $A$ is called nuclear
 if
$A\otimes_{\min} B=A\otimes_{\max} B$ for any
$C^*$-algebra $B$ (see [La]). For a long
time it remained an open problem whether it suffices for
the nuclearity
of $A$ to assume that
$$  A\otimes_{\min} A^{op}=A\otimes_{\max} A^{op}\leqno(*)$$
where $A^{op}$ is the opposite $C^*$-algebra (with the
product in reverse order).
In [Ki2], Kirchberg gave the first
counterexamples. However, he pointed out that it
remained unknown whether $(*)$ holds
in the (non-nuclear) case of $A=B(H)$. (Note that $B(H)$
is isomorphic to its opposite.)
Kirchberg also proposed an approach to this question
together with a series of equivalent conjectures.
One of his conjectures was our main motivation to
investigate the non-separability of the metric space
$OS_n$, and as a result, we  obtain a negative answer
to the above mentioned question: the identity
$(*)$ fails for $A=B(H)$.
In
other words, we have
 $$B(H)\otimes_{\min}B(H)\not=B(H)\otimes_{\max}B(H).$$
Equivalently, there is more than
one $C^*$-norm on
$B(H)\otimes B(H)$ whenever $H$ is an
 infinite dimensional Hilbert space.
This is proved in section 3.

In the same section, we also show that if $E$ is an
$n$-dimensional Banach space
then the operator space $\max(E)$ in the sense of
[BP] satisfies
$$d_{SK}(\max(E))\geq{\sqrt{n}\over 4}.$$
In particular, for any infinite dimensional Banach space $X$
the operator space
$\max(X)$ is not exact.

 After circulating a  first version of this paper, we
observed that the non-separability of $OS_n$ for  $n>2$
can be alternately deduced from properties of
Kazhdan-groups, following the ideas of Voiculescu [V] (see
Remark 2.10 below). However, this approach does not seem
to give quite as much information as our original one.

In the final section 4,
we introduce a "new"
tensor product obtained as follows.
Let $E\subset B(H)$ and $F\subset B(K)$ be two
 operator spaces,
we denote by $E\otimes_M F$ the completion of
the linear  tensor product with respect to the norm
induced on it by $B(H)\otimes_{\max} B(K)$. This
  tensor product makes sense both in the category
of operator spaces and in that of $C^*$-algebras.
Our previous results show that it differs
in general from the minimal tensor product. We include
several properties of this tensor product,
based mainly on [Ki3].

In the rest of this introduction we give
 some background and explain our notation.
We refer to [Sa, Ta] for operator algebra theory.
\n  Let $H,K$ be Hilbert spaces.
We will denote by $H\otimes_2 K$ their Hilbertian
 tensor product. By an
operator space we mean a closed subspace of $B(H)$ for
some Hilbert space $H$.
If $E_1\subset
 B(H_1)$, $E_2\subset B(H_2)$ are    operator spaces, we
will denote   by $E_1\otimes_{\min} E_2$
 their minimal (or spatial) tensor product equipped with
the minimal (or spatial) tensor norm
 induced by the space $B(H_1\otimes_2 H_2)$. When $E_1$ and
$E_2$ are $C^*$-subalgebras then this norm is a $C^*$-norm, and
actually
it is the smallest $C^*$-norm on the algebraic tensor
product $E_1\otimes  E_2$. In the case of $C^*$-algebras,
we will denote by $E_1\otimes_{\max} E_2$
the completion of $E_1\otimes  E_2$ with respect to the
largest $C^*$-norm on $E_1\otimes  E_2$.

 We recall that if $E\subset B(H)$ and $F\subset B(K)$ are operator
 spaces,
 then a map $u\colon
 \ E\to F$ is completely bounded (in short $c.b.$) if the
maps $u_m  =
 I_{M_m} \otimes u\colon \ M_m(E)\to M_m(F)$ are uniformly bounded when
 $m\to \infty$, i.e. if we have $\sup\limits_{m\ge 1}
\|u_m\|<\infty$. The
 $c.b.$ norm of $u$ is defined as
 $$\|u\|_{cb} = \sup_{m\ge 1} \|u_m\|.$$
\n We will say that $u$ is  completely isometric (resp.
completely contractive) or is a complete isometry (resp.
a complete contraction) if the maps $u_m$ are isometries
(resp. of norm $\le 1$) for all $m$. This is the same as
saying that for any operator space $X$ the map $$I_{X}\otimes u\colon
{X\otimes_{\min} E}\to X\otimes_{\min} F$$ is an
isometry (resp. a contraction). We also recall
 that $u$ is called completely
positive (in short $c.p.$) if all the maps $u_m$ are positive.

If $E_1,F_1$,
$E_2,F_2$ are operator spaces and if
$$u_1\colon\ E_1\to F_1\qquad u_2\colon\ E_2\to F_2$$
are completely bounded, then
$$u_1\otimes u_2\colon \ E_1\otimes_{\rm min} F_1 \longrightarrow
E_2\otimes_{\rm min} F_2$$
is $c.b.$\ and we have
$$\|u_1\otimes u_2\|_{cb} = \|u_1\|_{cb} \|u_2\|_{cb}.\leqno (0.2)$$

\n If $u\colon \ E_1\to F_1$ and
$u_2\colon\ E_2\to F_2$ are completely isometric then
$u_1\otimes u_2\colon \ E_1\otimes_{\rm min} E_2\longrightarrow
F_1\otimes_{\rm min} F_2$
also is completely isometric.
In particular, we note the completely isometric identity
$$M_n(E) = M_n \otimes_{\rm min} E.$$
We will use (0.2)
repeatedly in the sequel with no further reference.
We refer the reader to [Pa1] for more information.

It is known that the analogue of (0.2) fails for the $\max$-tensor
product. Instead we have:

\n (0.3) If $E_1,F_1$ and $E_2,F_2$ are
$C^*$-algebras (or merely operator systems) and if
$u_1\colon\ E_1\to F_1$ and $u_2\colon\ E_2\to F_2$
are completely positive maps, then the map
$u_1\otimes u_2$ defined on the algebraic tensor product
extends to a completely positive and  $c.b.$ map
$$u_1\otimes u_2\colon \ E_1\otimes_{\rm max} F_1 \longrightarrow
E_2\otimes_{\rm max} F_2$$
satisfying
$$\|u_1\otimes u_2\|_{cb(E_1\otimes_{\rm max} F_1 ,
E_2\otimes_{\rm max} F_2)} = \|u_1\| \|u_2\|.$$
See e.g. [Ta, Proposition 4.23],
 [Pa1, Proposition 3.5 and Proposition 10.11]
or [W2, Proposition 1.11] for more details.

We will use the duality theory for operator spaces,
which was introduced
in [BP, ER1] using Ruan's "abstract" characterization of
operator spaces ([Ru]). This can be summarized
as follows: Let $E$ be an operator space
and let $E^*$ be the dual Banach space. Then,
for some Hilbert space $K$, there is
an isometric embedding $E^*\subset B(K)$
such that for any operator space
$F$, the minimal (=spatial) norm
on
$E^*\otimes F$
coincides with the norm induced on it  by
the space $cb(E,F)$. Moreover, this property
characterizes the operator space
$$ E^*\subset B(K)$$
up to complete isometry.
We will refer to this operator space
as the ``operator space dual" of $E$.
We refer to [BP, ER1-2, B1-2, Ru] for detailed information.

Consider arbitrary Banach spaces $E,Y$ and an operator $u\colon \ E\to
Y$.

 Recall that $u$ is called 2-absolutely summing if there is a constant
 $C$
 such that for all $n$ for all $(x_1,\ldots, x_n)$ in $E$ we have
 $$\sum \|ux_i\|^2 \le C^2 \sup\left\{\sum |\xi(x_i)|^2\mid
\xi \in E^*,\ \
 \|\xi\|\le 1\right\}.$$
 We denote by $\pi_2(u)$ the smallest constant $C$ for which this
 holds.

 It is easy to check that for any bounded operators $v\colon \ Y\to
 Y_1$ and
 $w\colon\ E_1\to E$ we have
 $$\pi_2(vuw) \le \|v\|\pi_2(u)\|w\|.\leqno (0.4)$$

 \n {\bf Notation:}\ Let $(x_i)$ be a finite sequence in a
 $C^*$-algebra $A$
 we will denote for brevity
 $$RC((x_i)) = \max\left\{\left\|\sum x_ix_i^*\right\|, \left\|\sum
 x_i^*x_i\right\|\right\}.$$

Note that if $A$ is commutative and if $E\subset A$ is any subspace
 equipped with the induced norm, we have for all $(x_1,\ldots, x_n)$ in
 $E$
 $$RC((x_i)) = \sup\left\{\sum|\xi(x_i)|^2 \mid \xi\in E^*,
\ \ \|\xi\|\le
 1\right\}.\leqno (0.5)$$

 \proclaim Definition 0.3. Consider a $C^*$-algebra $A$ and a Banach
 space
 $Y$. Let $E\subset A$ be  subspace. A linear map $u\colon \ E\to Y$
 will be
 called $(2,RC)$-summing if there is a constant $C$ such that for all
 $n$
 and for all $(x_1,\ldots, x_n)$ in $E$ we have
 $$\sum\|ux_i\|^2 \le C^2 RC((x_i)).$$

 We denote by $\pi_{2,RC}(u)$ the smallest constant $C$
for which this holds.
 We refer to [P3] for a more systematic treatment of $(2,w)$-summing
 operators when $w$ is a ``weight'' in the sense of [P3]. (See also
 [P4, section 2] or
 [P2].) By a well known variant of a Pietsch's
factorization theorem for
  2-absolutely summing operators, we have

 \proclaim Proposition 0.4. Consider $u\colon \ E\to Y$ as in the
 preceding
 definition. The following are equivalent.\medskip
 \item{\rm (i)} $u$ is $(2,RC)$-summing and $\pi_{2,RC}(u)\le C$.
 \item{\rm (ii)} There are states $f,g$ on $A$ and $0\le \theta \le 1$
 such that
 $$\forall x\in E\qquad \|u(x)\| \le C\{\theta f(x^*x) + (1-\theta)
 g(xx^*)\}^{1/2}.$$
 \item{\rm (iii)} The map $u\colon\ E\to Y$ admits an extension $\tilde
 u\colon
 \ A\to Y$ such that $\pi_{2,RC}(\tilde u) \le C$.

 \pf This is -- by now -- a well known application of the Hahn-Banach
 theorem. For more details we refer the reader e.g.\ to
[P4, Lemma 1.3] or [P3, Prop. 5.1].

 For bilinear forms, we have the following known analogous
 statement.(In the commutative case,
this can be found in [Kw].)

 \proclaim Proposition 0.5. Let $A,B$ be $C^*$-algebras and let
 $E\subset
 A$, $F\subset B$ be closed subspaces. Let $C>0$ be a fixed constant.
 The
 following properties of a linear map $u\colon \ E\to F^*$ are
 equivalent.
 \medskip
 \item{\rm (i)} For any $n$, any $(x_1,\ldots, x_n)$ in $E$ and
any
 $(y_1,\ldots, y_n)$ in $F$, we have
 $$\left|\sum \langle u(x_i), y_i\rangle\right| \le C[RC((x_i))
 RC((y_i))]^{1/2}.$$
 \item{\rm (ii)} There are states $f_1,g_1$ on $A$, $f_2$, $g_2$ on $B$
 and
 $0\le \theta_1, \theta_2 \le 1$ such that
 $$\eqalign{\forall (x,y)\in E\times F\qquad |\langle u(x),y
 \rangle|&\le
 C[\theta_1f_1(x^*x) + (1-\theta_1)g_1(xx^*)]^{1/2} \cr
 &\qquad [\theta_2f_2(y^*y) +
 (1-\theta_2)g_2(yy^*)]^{1/2}.}$$
 \item{\rm (iii)} For some Hilbert space $H$, $u$ admits a
 factorization of the
 form
 $$E {\buildrel a\over \longrightarrow} H {\buildrel b^*\over
 \longrightarrow} F^*$$
 with operators $a\colon \ E\to H$ and $b\colon \ F\to H^*$ such that
 $$\pi_{2,RC}(b) \pi_{2,RC}(a) \le C.$$

 \pf See e.g.\ [P3, Theorem 6.1] or [P4, Lemma 1.3].

 \n {\bf Notation:}\ Let $u\colon\ X\to Y$ be an operator between
 Banach
 spaces. Assume that $u$ factors through a Hilbert space $H$, i.e.\ we
 have
 $u=\alpha\beta$ with $\alpha\colon \ H\to Y$ and $\beta\colon \ X\to
 H$.
 Then we will denote
 $$\gamma_2(u) = \inf\{\|\alpha\|\ \|\beta\|\}$$
 where the infimum runs over all possible factorization. This is the
 ``norm
 of factorization through Hilbert space'' of $u$. See [P2] for more on
 this theme.

 \proclaim Corollary 0.6. Let $u$ be as
 in Proposition 0.5. Then if (iii)
 holds there
are operators
 $\tilde a\colon \ A\to H$ and $\tilde b\colon \ B\to H^*$
such that  $\tilde b^*\tilde a$ (viewed as a bilinear form on $A\times
B$)
extends $u$ and
 $$\pi_{2,RC}(\tilde a) \pi_{2,RC}(\tilde b) \le C.$$
In particular, if we let $\tilde u= \tilde b^*\tilde a$, then
the operator $\tilde u\colon \ A\to B^*$
satisfies $\|\tilde u\|
 \le C$ and $\langle\tilde u(x),y\rangle = \langle u(x),y\rangle$ for
 all $(x,y)$ in $E\times F$.  Moreover we have $\gamma_2(u)\le C$.

\pf This follows from Propositions 0.4 and 0.5.

The following fact is well known to specialists. (See e.g.\ [Ki1] Lemma
3.9)

\proclaim Lemma 0.7. Let $C$ be a $C^*$-algebra and let $I\subset C$ be
a
closed ideal. Let $E,X$ be operator spaces with $E\subset X$. Consider
the
canonical (complete) contraction $C\otimes_{\rm min} X\to
(C/I)\otimes_{\rm
min} X$. Since this map vanishes on $I\otimes_{\rm min} X$, we clearly
have a
canonical (complete) contraction
$$T_X\colon \ C\otimes_{\rm min} X/I\otimes_{\rm min} X\to (C/I)
\otimes_{\rm min} X.$$
If $T_X$ is an isomorphism, then $T_E$ also is an isomorphism,
moreover
$$\|T^{-1}_E\| \le \|T^{-1}_X\|.$$

\pf It is well known (\cf [Ta, p. 27]) that $I$ possesses
 an approximate unit formed of
elements $p_i\in I$ such that $0\le p_i\le 1$ (hence $\|1-p_i\|\le 1$)
and
$p_ix\to x$ $\forall x\in I$. Then the proof can be completed exactly
as in
[P1, Lemma 3].

\n {\bf Remark.} Equivalently, if we consider the complete isometry
$$C\otimes_{\rm min} E\to C\otimes_{\rm min} X$$
then this map defines after passing to the quotient spaces a complete
isometry
$$C\otimes_{\rm min} E/I\otimes_{\rm min} E\to C\otimes_{\rm min} X/I
\otimes_{\rm min} X.$$

We will also invoke the following elementary fact
which follows easily (like the preceding result) from
the existence of an approximate unit in any ideal
of a $C^*$-algebra.
(Recall that $A\otimes B$ denotes the algebraic tensor
product.)
\proclaim Lemma 0.8. Let $A,B,C$ be $C^*$-algebras.
Let $\pi\colon C\to A$ be a surjective $*$-homomorphism. Let
${\cal I}=Ker(\pi)$.
 Then (viewing the three sets appearing below as subsets
of $C\otimes_{\min} B$) we have
$$[C\otimes B ]\cap [{{\cal I}}\otimes_{\min} B]={{\cal I}}\otimes B.$$
Equivalently, let
$$q\colon C\otimes_{\min} B \to [C\otimes_{\min} B]/
[{\cal I}\otimes_{\min} B]$$
be the quotient map and let
$$T\colon [C\otimes_{\min} B]/
[{\cal I}\otimes_{\min} B] \to A\otimes_{\min} B$$
be the morphism associated to $\pi\otimes I_B$.
Then $T$ induces a linear and $*$-algebraic isomorphism
between $q(C\otimes B)$ and $A\otimes B$.

\n {\bf Acknowledgement:} The second author
 would like to thank  E. Kirchberg
for introducing him to the questions
considered in this paper and
B. Maurey for a
stimulating conversation.

\vfill\eject

\n {\bf \S 1. Factorization of bilinear forms on exact operator
spaces.}

Let $n\ge 1$ be an integer. We will denote by
$$J_n\colon \ M_n\to M^*_n$$
the map defined by
$$\forall x,y\in M_n\quad \langle J_n(x),
 y\rangle = {1\over n} tr(^tyx).$$
The following notion is natural for our subsequent results. It has
already
been considered in [I].

\proclaim Definition 1.1. Let $E,F$ be operator spaces.
Let $u\colon\ E\to
F^*$ be a linear map. We will say that $u$ is tracially bounded (in
short
t.b.) if
$$\sup_{n\ge 1} \|J_n\otimes u\|_{M_n(E)\to M_n(F)^*} < \infty$$
and we denote
$$\|u\|_{tb}  = \sup_{n\ge 1}\|J_n\otimes u\|_{M_n(E)\to M_n(F)^*}.$$
Equivalently, $u\colon \ E\to F^*$ is t.b.\ iff the bilinear forms
$$u_n\colon \ M_n(E) \times M_n(F)\to {\bf C}$$
defined by
$$u_n((x_{ij}), (y_{ij})) = {1\over n} \sum_{ij} \langle u(x_{ij}),
y_{ij}\rangle$$
are bounded uniformly in $n$ and we have
$$\|u\|_{tb} = \sup_{n\ge 1}\|u_n\|.$$

We immediately observe

\proclaim Lemma 1.2. For a linear map $u\colon \ E\to F^*$
$${\it complete\
boundedness} \Rightarrow {\it tracial\  boundedness} \Rightarrow
{\it boundedness},$$ and we have
$$\|u\| \le \|u\|_{tb} \le\|u\|_{cb}.$$

\pf If $u$ is $c.b.$\ then for any $u$ and any $(x_{ij})\in M_n(E)$
$$\|(u(x_{ij}))\|_{M_n(F^*)} \le \|u\|_{cb} \|x\|_{M_n(E)}.$$
We have $M_n(F^*) = cb(F,M_n)$, hence
$$\eqalign{\|(u(x_{ij}))\|_{M_n(F^*)} &= \sup\{\|\langle u(x_{ij}),
y_{k\ell}\rangle\|_{M_n(M_n)}\ | \  \|(y_{k\ell})\|_{M_n(F)} \le 1\}\cr
&= \sup_{y\in B_{M_n(F)}}
\left\{\left|\sum_{ijk\ell} \langle u(x_{ij}),
y_{k\ell}\rangle \alpha_{j\ell} \beta_{ik}\right|\ \Big|\
\sum_{j\ell} |\alpha_{j\ell}|^2 \le 1\ \sum_{ik}
|\beta_{ik}|^2\le 1\right\}.}$$
Taking $\alpha,\beta = {I\over \sqrt n}$ we
get
$$\|(u(x_{ij}))\|_{M_n(F^*)} \ge \sup\left\{{1\over n}\left|\sum
\langle
u(x_{ij}), y_{ij}\rangle\right|\ \Big|\ \|(y_{ij})\| _{M_n(F)}\le
1\right\}
$$
hence $\|u\|_{tb}\le \|u\|_{cb}$. The inequality $\|u\| \le \|u\|_{tb}$
is
clear by taking $n=1$.\qed

The following consequence of the
 non-commutative Grothendieck inequality
is known ([B3]).

\proclaim Lemma 1.3. Let $A,B$ be $C^*$-algebras. Then any bounded
linear
operator $u\colon\ A\to B^*$ is tracially bounded and we have
$$\|u\|\le \|u\|_{tb} \le K\|u\|$$
for some numerical constant $K$. Let $E\subset A$ and $F\subset B$ be
closed subspaces and let $i_E\colon \ E\to A$, $i_F\colon \ F\to B$ be
the
inclusion maps. Then the restriction $i^*_Fui_E\colon \ E\to F^*$
satisfies
$$\|i^*_Fui_E\|_{tb} \le K\|u\|.$$

\pf Consider $x=(x_{ij}) \in M_n(A)$ and $y = (y_{ij}) \in M_n(B)$.

Let $\alpha(x) = \max\left\{n^{-1/2}\left\|\sum\limits_{ij}
x^*_{ij}x_{ij}
\right\|^{1/2}, n^{-1/2}\left\|\sum\limits_{ij} x_{ij}x^*_{ij}\right\|
^{1/2}\right\}$. It is easy to check that
$$\alpha(x) \le \|x\|_{M_n(A)}.$$
By the non-commutative Grothendieck inequality (cf.\ [H1], see also
[P2])
we have for some numerical absolute constant $K$
$$\eqalignno{\left|{1\over n} \sum_{ij}\langle u(x_{ij}), y_{ij}\rangle
\right| &\le K\|u\| \alpha(x) \alpha(y)\cr
\noalign{\hbox{hence a fortiori}}
&\le K\|u\| \|x\|_{M_n(A)} \|y\|_{M_n(B)}.}$$
Therefore $\|u\|_{tb} \le K\|u\|$. The second part is obvious since
$$\|i^*_Fui_E\|_{tb} \le \|u\|_{tb}.$$\qed

\n {\bf Remark.} The preceding argument actually shows
the following: Consider
$u\colon A\to B^*$ of the form $u=a^*b$, where $b\colon
E\to H$ and $a\colon F\to H^*$ are $(2,RC)-$summing
operators. Then we have
$$\|u\|_{tb} \le \pi_{2,RC}(a)\pi_{2,RC}(b).$$

We will now show that if $E$ and $F$ are exact operator spaces, the
converse to the second part of Lemma~1.3 also holds,
that is to say the bilinear form on $E\times F$ associated to
a tracially bounded map
$u\colon E\to F^*$ is the restriction
of a bounded bilinear form on $A\times B$.
 This is the key result
for this paper.

\n We denote by $F_\infty$ the free group with countably many
generators
denoted by $g_1,g_2,\ldots \ .$ Let $\lambda\colon \ F_\infty \to
B(\ell_2(F_\infty))$ be the left regular representation and let us
denote
simply by $C_\lambda$ the reduced $C^*$-algebra of $F_\infty$. Let $E$
be
an operator space. Let $(x_t)_{t\in F_\infty} $ be a
finitely supported family of elements of $E$. For simplicity we will
denote simply by
$$\left\|\sum \lambda(t) \otimes x_t\right\|_{\rm min}$$
the norm
 of $\sum_{t\in F_\infty} \lambda(t)\otimes x_t$
in $C_\lambda\otimes_{\rm
min} E$.

The following inequality (cf.[HP, Prop.~1.1])
plays an important r\^ole in the sequel.
 $$\hbox{For any finite sequence $(x_i)$ in a $C^*$-algebra we have}
 \leqno
 (1.1)$$
 $$\left\|\sum \lambda(g_i)\otimes x_i\right\|_{\rm min} \le 2
 \max\left\{
 \left\|\sum x^*_ix_i\right\|^{1/2}, \left\|\sum x_ix^*_i\right\|^{1/2}
 \right\}.$$

\proclaim Theorem 1.4. Let $E,F$ be exact operator spaces. Let $C =
d_{SK}(E) d_{SK}(F)$. Let \break $u\colon \ E\to F^*$ be a tracially
bounded
linear map. Let $(x_t)_{t\in F_\infty}$
 (resp. $(y_t)_{t\in F_\infty}$) be a finitely supported
family of elements of $E$ (resp. $F$). Then we have
$$\left|\sum \langle u(x_t), y_t\rangle \right|
\le C\|u\|_{tb} \left\| \sum \lambda(t) \otimes
x_t\right\|_{\rm min} \left\|\sum \lambda(t) \otimes
y_t\right\|_{\rm min},\leqno
(1.2)$$
In particular,  for all $n$ and all $x_i \in E$, $y_i \in F$
($i=1,2,...,n$) we have
$$\left|\sum \langle u(x_i), y_i\rangle \right|
\le C\|u\|_{tb} \left\| \sum \lambda(g_i) \otimes
x_i\right\|_{\rm min} \left\|\sum \lambda(g_i) \otimes
y_i\right\|_{\rm min},\leqno
(1.3)$$
hence
$$
\le4C\|u\|_{tb}\max\left\{\|\sum
x_i^*x_i\|^{1/2},\|\sum x_ix_i^*\|^{1/2}\right\}\max\left\{\|\sum
y_i^*y_i\|^{1/2},\|\sum y_iy_i^*\|^{1/2}\right\}. \leqno
(1.3)'$$
Furthermore, if $A,B$ are $C^*$-algebras  such that $E\subset A$,
$F\subset B$ (completely isometrically) then $u$ admits an extension
$\tilde u\colon\ A\to B^*$ such that
$$\|\tilde u\|\le  \|\tilde u\|_{tb}\le 4C\|u\|_{tb}$$
and $\langle \tilde u(x), y\rangle = \langle u(x), y\rangle$ $\forall
(x,y)\in E\times F$.

\pf The proof is essentially the same as that of in [P1, Theorem 8]. An
essential
ingredient in the proof is Wassermann's construction of a specific
embedding of $C_\lambda$ into an ultrapower (in the von~Neumann sense)
of
matrix algebras, cf.\ [W1]. This is based on the residual finiteness
of the free group.

\n More precisely, consider  the family
$$\{M_\alpha\mid \alpha \ge 1\}$$
of all matrix algebras. Let $L = \{(x_\alpha)_{\alpha\ge 1}\mid
x_\alpha
\in M_\alpha,\  \sup\limits_{\alpha\ge 1}\|x_\alpha\|_{M_\alpha} <
\infty\}$
equipped with the norm $\|(x_\alpha)\|_L = \sup\limits_{\alpha \ge 1}
\|x_\alpha\|_{M_\alpha}$. Let ${\cal U}$ be an ultrafilter on $\bf N$.
Let
us denote by $\tau_\alpha$ the normalized trace on $M_\alpha$, and let
$$I_{\cal U} = \{(x_\alpha)_{\alpha\ge 1}\in L\mid \lim_{\cal U}
\tau_\alpha(x^*_\alpha x_\alpha) = 0\}.$$
We then set ${\cal N} = L/I_{\cal U}$. It is a well known fact
(cf\ [Sa]) that ${\cal N}$
is a finite von~Neumann algebra. Let $VN_\lambda$ be the von~Neumann
algebra
generated by $\lambda$. (Note that $C_\lambda\subset VN_\lambda$.)
Wassermann proved that for a suitable ${\cal U}$ one can
 find for each
$g$ in $F_\infty$ a sequence $(u^g_\alpha)_{\alpha\ge 1}$ such
that:\medskip

\item{\rm (i)} $u^g_\alpha$ is unitary in $M_\alpha$ and has real
entries $(\alpha\ge 1, i=1,2,\ldots)$,
\item{\rm (ii)} $\lim\limits_{\cal U} \tau_\alpha({u^s_\alpha}^{*}
u^g_\alpha) = 0$
if $g\ne s$, or equivalently since the entries are real
$\lim\limits_{\cal U} \tau_\alpha(^t{u^s_\alpha} u^g_\alpha) = 0$ if
$g\ne s$.
\item{\rm (iii)} The mapping $\Phi\colon \ VN_\lambda \to L/I_{\cal U}$
which
takes $\lambda(g)$ to the equivalence class of $(u^g_\alpha)_{\alpha\ge
1}$ is an isometric representation mapping $VN_\lambda$ onto a
von~Neumann
subalgebra of\break ${\cal N} = L/I_{\cal U}$. A fortiori it is
completely isometric.

The last point implies that we can write
$$\left\|\sum_{t \in  F_\infty} \lambda(t) \otimes x_t\right\|_{\rm
min} = \left\|
\sum_{t \in  F_\infty} \Phi[\lambda(t)]\otimes x_t\right\|_{(L/I_{\cal
U})\otimes _{\rm
min}E}.\leqno (1.4)$$
 Without loss
of generality we may assume $\dim E = \dim F = n$. Hence, for any
$\vp>0$, for some
integer  $N$ there
is $E_1\subset M_N$ such that
$$d_{cb}(E, E_1) < d_{SK}(E) (1+\varepsilon).\leqno (1.5)$$
Clearly we have completely isometrically
$$L/I_{\cal U} \otimes_{\rm min} M_N = (L\otimes_{\rm min} M_N)/I_{\cal
U}
\otimes_{\rm \min} M_N.$$
By Lemma 0.7, since $I_{\cal U}$ is an ideal this remains true with
$E_1$
in the place of $M_N$. By (1.5), it follows that the natural
(norm one)  map
$$T_E\colon \ (L\otimes_{\rm min} E)/(I_{\cal U}\otimes_{\rm min} E)
\to
(L/I_{\cal U}) \otimes_{\rm min} E$$
has an inverse with norm
$\|T^{-1}_E\| < d_{SK}(E)(1+\varepsilon)$, and since
$\varepsilon$ is arbitrary and the same is true for $F$,
 we actually have
$$\|T^{-1}_E\| \le d_{SK}(E)\qquad \|T^{-1}_F\| \le
d_{SK}(F).\leqno(1.6)$$
Since $u$ is tracially bounded, for any $\alpha$, the linear map
\def\a{\alpha}
$V_\alpha\colon M_\a(E)\to M_\a(F)^*$ defined by
$$\langle V_\a(a_\a\otimes x), b_\a\otimes y\rangle =\tau_\a(^t b_\a
a_\a) \langle u(x),y\rangle ,\quad
\forall x\in E,\ \forall y\in F,\ \forall a_\a,\  b_\a\in M_\a
$$
is bounded with $$\|V_\a\|\le \|u\|_{tb}.\leqno(1.7)$$
For $a\in L$, let us denote by $(a_\a)$ its coordinates.
Similarly consider an element
$z$ in $L\otimes E$. Clearly $z$ can be identified to a family
$(z_\a)$ with $z_\a\in M_\a\otimes E$. Note that
$$\|z\|_{L\otimes_{\min} E} = \sup_{\a}
\|z_\a\|_{M_\a(E)}.\leqno(1.8)$$
With this notation, we can define a linear  map
$V\colon L\otimes_{\min} E \to (L\otimes_{\min} F)^*$
by setting
$$\langle V(a\otimes x), b\otimes y\rangle =\lim_{\cal U}
\tau_\a(^tb_\a a_\a) \langle u(x),y\rangle \qquad
\forall a,b \in L\ \forall x\in E,\ \forall y\in F
.$$
Clearly, by (1.7) and (1.8)  we have
$\|V\|\le \sup_\a\|V_\a\| \le \|u\|_{tb}.$ Moreover, it is clear that
for all $\xi$ in $I_{\cal U}\otimes_{\rm min} E$ and all $\eta$ in
$I_{\cal U}\otimes_{\rm min} F$
we have $\langle V(\xi),\eta\rangle =0$.
Therefore $V$ defines canonically a map
$$\tilde V\colon (L\otimes_{\rm min} E)/(I_{\cal U}\otimes_{\rm min}
E)  \to [(L\otimes_{\rm min} F)/(I_{\cal U}\otimes_{\rm min} F)]^*$$
such that $\|\tilde V\| = \|V\|\le \|u\|_{tb}$.
By (1.6), $\tilde V$ also defines a map
$$\hat V\colon (L/I_{\cal U}) \otimes_{\rm min} E\to
[(L/I_{\cal U}) \otimes_{\rm min} F]^* $$
such that $$\|\hat V\|\le  C\|u\|_{tb}.\leqno(1.9)$$
Let $(x_t)$ and $(y_t)$ be as in Theorem 1.4, and
 let
$T_1=\sum \lambda(t)\otimes x_t$,\qquad
$T_2=\sum \lambda(t)\otimes y_t$, and let
$$\hat{T_1}=(\Phi\otimes I_E) (T_1) \in   (L/I_{\cal U}) \otimes_{\rm
min} E,
\
\hat{T_2}=(\Phi\otimes I_F) (T_2) \in   (L/I_{\cal U}) \otimes_{\rm
min} F.$$
 By (ii) above, we clearly have
$\sum \langle u(x_t),y_t\rangle = \langle \hat V (\hat{T_1}),
\hat{T_2}\rangle $
hence
$$|\sum \langle u(x_t),y_t\rangle |\le \|\hat V\| \|\hat{T_1} \|_{\min}
\|\hat{T_2} \|_{\min},$$
and this together with (1.4) and  (1.9)
implies
(1.2). This proves (1.2). Clearly (1.3) is but a
particular case of (1.2) and (1.3)' follows, using (1.1).

\n Finally, the last assertion follows
from Corollary 0.6 and the remark after Lemma 1.3.\qed

\n {\bf Remark.} The preceding proof of (1.2) remains
valid if we replace $F_\infty$ by any  residually
finite discrete group $G$. Actually, we only use the
fact that there is a  completely isometric embedding into
an ultraproduct, say
$\Phi\colon C^*_\lambda(G)\to L/I_{\cal U}$
such that
$$\forall s,g\in G\quad \lim_{\cal U}\tau_\alpha
(^t\Phi_\a(\lambda(s))\Phi_\a(\lambda(g)))=\delta_{sg}.$$

A corollary of the non-commutative Grothendieck theorem says that every
 bounded linear operator $u\colon \ A\to B^*$ factors through a Hilbert
 space
 when $A,B$ are $C^*$-algebras and we have $\gamma_2(u) \le K\|u\|$ for
 some
 absolute constant $K$. In the same vein, we have

 \proclaim Corollary 1.5. Let  $E,F$ be as in Theorem~1.4 with
 $C=d_{SK}(E)
 d_{SK}(F)$. Then, \break if $u\colon \ E\to F^*$ is
tracially bounded  there
 is a Hilbert space $H$ and a factorization $u=a^*b$
 $$E {\buildrel b\over \longrightarrow} H {\buildrel a^*\over
 \longrightarrow} F^*$$
 with $a,b$ such that
 $$\pi_{2,RC}(a) \pi_{2,RC}(b) \le 4C\|u\|_{tb}.$$
 Conversely, if there is such a factorization we have
 $$\|u\|_{tb} \le \pi_{2,RC}(a) \pi_{2,RC}(b).\leqno (1.10)$$

 \pf The direct implication follows from Theorem~1.4 and
 Proposition~0.5. The
 converse follows from the remark after Lemma~1.3.\qed

 \proclaim Corollary 1.6. Let $E,F$ be exact operator spaces. Let
 $C=d_{SK}(E) d_{SK}(F)$. Then every completely bounded map $u\colon
 \ E\to
 F^*$ factors through a Hilbert space and we have
 $$\gamma_2(u) \le 4C\|u\|_{cb}.\leqno (1.11)$$

 \pf First recall that $\|u\|_{tb} \le \|u\|_{cb}$ by Lemma~1.2. Then
 this
 is deduced from Theorem~1.4 using Proposition~0.5 and
 Corollary~0.6.\qed

 \proclaim Corollary 1.7. Let $E$ be a closed subspace of a commutative
 $C^*$-algebra and let $F$ be an exact operator space.
Then every completely bounded
 $u\colon \ E\to F^*$ is 2-absolutely summing and we have $\pi_2(u) \le
 4d_{SK}(F)\|u\|$.

 \pf Since a commutative $C^*$-algebra is nuclear we have $d_{SK}(E)
 =1$.
 The result then follows from Theorem~1.4 taking (0.5) into
 account.\qed

 In particular we have the following corollary which is already known.
 It
 was proved independently by V.~Paulsen and the second author on one
 hand
  (using Clifford matrices, this gives the better constant 2, see
  [Pa3])
 and by the first author on the other
(using random matrices, this yields a worse constant).

 \proclaim Corollary 1.8. Consider $C^*$-algebras $A,B$ and subspaces
 $E\subset A$ and $F\subset B$. If $A,B$ are assumed commutative, then
 any
 $c.b.$ map $u\colon\  E\to F^*$ can be written as $u=a^*b$ with
 $\pi_2(a)
 \pi_2(b) \le 4\|u\|_{cb}$.

 \pf Again, by the nuclearity of $A,B$ we have $d_{SK}(E) = d_{SK}(F)
 =1$.
 Note that by (0.5) we have $\pi_2(a) = \pi_{2,RC}(a)$ for
all $a\colon \
 E\to H$ and similarly for $F$. Hence this follows from Theorem~1.4.
 \qed

 \proclaim Corollary 1.9. Let $E$ be an exact operator space. Consider
 a
 linear map $v\colon \ E\to H$ into a Hilbert space. Let $C>0$ be a
 constant.
 Assume that for all $n$ and all $(x_{ij})$ in $M_n(E)$ we have
 $$\left({1\over n} \sum_{ij} \|v(x_{ij})\|^2\right)^{1/2}
 \le C \|(x_{ij})\|_{M_n(E)}.
 \leqno (1.12)$$
 Then $v$ is $(2,RC)$-summing and
 $$\pi_{2,RC}(v) \le 2Cd_{SK}(E).\leqno (1.13)$$

 \pf Consider the mapping $u = \overline v^*v\colon\ E\to \overline
 E^*$,
 obtained by identifying $H$ with its antidual $\overline H^*$. From
 (1.12) it is easy to deduce by Cauchy-Schwarz that $\|\overline
 v^*v\|_{tb}
 \le C^2$. Hence by (1.2) we have for all $(x_i)$ in $E$
 $$\sum \|vx_i\|^2 = \sum \langle \overline v^*vx_i,x_i\rangle \le 4C^2
 d_{SK}(E)^2 RC((x_i))^2$$
 and (1.13) follows.\qed

\def\p.{\medskip{\noindent{\bf Proof. }}}
We need to recall some elementary facts on ultraproducts of operator
spaces.
Let $(E_i)_{i\in I}$ be a family of operator spaces with $E_i\subset
B(H_i)$.
Then their ultraproduct $\hat E=\Pi E_i/{\cal U}$ embeds isometrically
into $\Pi B(H_i)/{\cal U}$. The latter being a $C^*$-algebra, this
embedding
defines an operator space structure on $\hat E$. It is easy to check
that we
have isometrically
$$M_n(\hat E)=\Pi M_n(E_i)/{\cal U}.$$
Equivalently
$$M_n\otimes_{\min}\hat E=\Pi (M_n\otimes_{\min}E_i)/{\cal U}.$$
This identity clearly remains valid with $M_n$ replaced by any subspace
$F\subset M_n$, therefore we also have the following isometric
identity, valid
if $d_{SK}(F)=1\  ( \dim F<\infty)$.
$$F\otimes_{\min}\hat E=\Pi (F\otimes_{\min}E_i)/{\cal U}.\leqno
(1.14)$$
This yields

\proclaim Corollary 1.10. Let $I$ by any set. Let $(E_i)_{i\in
I},\ (F_i)_{i\in
I}$ be exact operator spaces with
$$C=\sup_{i\in I}d_{SK}(E_i)d_{SK}(F_i)<\infty.$$
Let $u_i\colon E_i\to F_i^*$ be tracially bounded maps with $\sup_{i\in
I}\N{u_i}_{tb}\leq 1$.
Let ${\cal U}$ be an ultrafilter on $I$ and let $\hat E$ (resp. $\hat
F)$ be the
ultraproduct of $(E_i)_{i\in I}$ (resp. $(F_i)_{i\in I}$) modulo ${\cal
U}$. Let
$\hat u\colon \hat E\to(\hat F)^*$ be the map associated to the family
$(u_i)_{i\in I}$.
Then for all finite sets $(x_k)$ in $\hat E$ and $(y_k)$ in $\hat F$ we
have
$$|\sum\langle \hat u(x_k),y_k\rangle |\leq C\N{\sum\lambda(g_k)\otimes
x_k}_{\min}\N{\sum\lambda(g_k)\otimes y_k}_{\min}.\leqno (1.15)$$

\p. We use the fact that $C_\lambda^*(F_\infty)$ is an exact
$C^*$-algebra since
it has the slice map property (\cf [DCH, Corollary 3.12] and see [Kr]
for more
details). Therefore if $F={\rm span}
(\lambda(g_1),...,\lambda(g_n))$ we have $d_{SK}(F)=1$. Then it is easy
to
derive (1.15) from (1.2) (applied to each $u_i$) taking (1.14) into
account.\qed

\vfill\eject

\n {\bf \S 2. The non-separability of $OS_n$.}

We will denote by $OS_n$ the set of all $n$ dimensional operator
spaces. We
identify two elements $E,F\in OS_n$ if they are completely isometric.
For
$E,F\in OS_n$, let
$$d_{cb}(E,F) = \inf\{\|u\|_{cb} \|u^{-1}\|_{cb}\mid u\colon \ E\to F,
u
\hbox{ complete isomorphism}\}.$$
Then it can be shown (see [P1] for the easy details) that
$$\delta_{cb}(E,F) = \hbox{Log } d_{cb}(E,F)$$
is a distance on $OS_n$ for which it is a complete metric space.

We will need a weaker metric structure on the space $OS_n$. To
introduce it
we need the following notation:\ For any linear map $u\colon \ E\to F$
between operator spaces we denote
$$\|u\|_k = \|I_{M_k} \otimes u\|_{M_k(E)\to M_k(F)}.$$
Note that
$$\|u\|_{cb} = \sup_{k\ge 1}\|u\|_k.$$
Now consider $E,F\in OS_n$. We define
$$d_k(E,F)  = \inf\{\|u\|_k \|u^{-1}\|_k\mid u\colon \ E\to F, u \hbox{
linear isomorphism}\}.$$
Then by a simple compactness argument (the unit ball of the space
 ${\cal L}(E,F)$ of all linear maps is
compact for any norm on ${\cal L}(E,F)$) one can check that
$$d_{cb}(E,F) = \sup_{k\ge 1} d_k(E,F).\leqno (2.1)$$
We set
$$\delta_w(E,F) = \sum_{k\ge 1} 2^{-k} \hbox{ Log } d_k(E,F).$$
Then, $\delta_w$ is a distance on $OS_n$. Let $\{E_i\}$ be a sequence
in
$OS_n$. Then $\delta_w(E_i,E)\to 0$ iff
$$d_k(E_i,F)\to 1\quad \hbox{for all}\quad k\ge 1.$$
In that case we will write simply $E_i {\buildrel w\over
\longrightarrow}
E$. It was observed in [P1] that $E_i {\buildrel w\over
\longrightarrow}
E$ iff for any non-trivial ultrafilter ${\cal U}$ on $\bf N$ the
ultraproduct $\Pi
E_i/{\cal U}$ is completely isometric to $E$.

We will need the following known fact.

\proclaim Proposition 2.1. For any $E,F$ in $OS_n$ and any $k\ge 1$ we
have
$$\eqalignno{d_{cb}(E,F)  &= d_{cb}(E^*,F^*)\cr
\noalign{\hbox{and}}
d_k(E,F) &= d_k(E^*,F^*)}$$
for all $k\ge 1$. Hence in particular
$$\delta_{cb}(E,F) = \delta_{cb}(E^*,F^*)\quad \hbox{and}\quad
\delta_w(E,F) = \delta_w(E^*,F^*).$$

\pf It clearly suffices to know that for any $u\colon \ E\to F$ we have
$$\|u\|_{cb} = \|u^*\|_{cb}\leqno (2.2)$$
and
$$\forall k\ge 1\qquad \|u\|_k =\|u^*\|_k.\leqno (2.3)$$
The identity (2.2) was proved in [BP, ER1], while (2.3) is easy to
check
using the definition of $\|u\|_k$. We have (by [Sm])
$$\|u\|_k = \sup\{\|bua\|_{cb}\mid a\colon \ M^*_k \to E, b\colon\ F\to
M_k, \|a\|_{cb}\le 1, \|b\|_{cb}\le 1\}.\leqno (2.4)$$
Clearly (2.4) implies (2.3). Then the above proposition is obvious.\qed

The following was proved in [P1].

\proclaim Proposition 2.2. Let $E\in OS_n$. The following are
equivalent
\medskip
\item{\rm (i)} For any sequence $\{E_i\}$ in $OS_n$ tending weakly to
$E$ we
have $d_{cb}(E,E_i)\to 1$ when $i\to\infty$.
\item{\rm (ii)} Same as (i) with each $E_i$ isometric to $E$.
\item{\rm (iii)} $d_{SK}(E) = d_{SK}(E^*) = 1$.\medskip

\n {\bf Remark.} For any fixed integer $k\ge 1$, these are also
equivalent
to the same property as (i) restricted to $E_i$ $k$-isometric to $E$.

\pf We only prove (i) $\Rightarrow$ (iii) which is what we use in the
sequel.

\n Assume (i). Then $E\subset B(\ell_2)$. Let $P_m\colon\ B(\ell_2)\to
M_m$ be
the projection which maps $e_{ij}$ to itself if $1\le i,j\le m$ and to
zero
otherwise. Let $E_m = P_m(E)\subset M_m$. It is very easy to check that
$E_m {\buildrel w\over \longrightarrow} E$. Hence if (i) holds we have
$$\eqalign{d_{SK}(E) &\le d_{cb}(E,E_m) d_{SK}(E_m)\cr
&\le d_{cb}(E,E_m)\cr
\hbox{hence}\qquad\qquad d_{SK}(E) &\le \lim d_{cb}(E,E_m) =1.}$$
Now if we apply the same to $E^*$ (equipped with the dual operator
space
structure) we obtain by Proposition~2.1 that $d_{SK}(E^*)=1$.

\n  Reformulated
in more concise terms, the proof reduces to this:\ let $OS_n(m)$ be the
subset of $OS_n$ formed of all $n$-dimensional subspaces of $M_m$. Then
the
union $\bigcup\limits_{m\ge  n} OS_n(m)$ is weakly dense in $OS_n$.
Hence
if (i) holds, $E$ (and also
$E^*$  by Proposition~2.1)  must be in the strong closure of
$\bigcup\limits_{m\ge n}
OS_n(m)$, which means that $d_{SK}(E)=1$ (and  $d_{SK}(E^*)=1$).

 We can now prove

\proclaim Theorem 2.3. The metric space $(OS_n, \delta_{cb})$ is
non-separable if $n>2$.

\pf Let $f\colon \ (OS_n,\delta_w)\longrightarrow (OS_n, \delta_{cb})$
be
the identity mapping. Note that $f^{-1}$ is continuous, but in general
$f$
is not.

\n However, if we assume $OS_n$ strongly separable then we claim that
$f$ is
in the first Baire class. Indeed, by (2.1) for any closed ball $\beta$
in
$(OS_n, \delta_{cb})$, $f^{-1}(\beta)$ is weakly closed. Hence if
$OS_n$ is
strongly separable, for any $U$ strongly open in $OS_n$ $f^{-1}(U)$
must be
an $F_\sigma$-set in weak topology, hence $f$ is in the first Baire
class.
Note that the domain of $f$ is compact hence is a Baire space. By
Baire's
classical theorem (cf.\ [Ba], see also [Ku, 31, X, Th. 1, p. 394]), if
the range of $f$ is separable, the set of
points of continuity of $f$ must be dense in the domain of $f$,
i.e.\ dense
for the weak topology.

\n This implies by Proposition~2.2 that for any $E$ in $OS_n$ there is
a
sequence $\{E_i\}$ in $OS_n$ such that $d_{SK}(E_i) = d_{SK}(E^*_i)=1$
which tends weakly to $E$. Equivalently, $E$ can be viewed as
the ultraproduct of $(E_i)$ with respect to a non-trivial
ultrafilter ${\cal U}$.

\n By Corollary 1.10  applied to the identity of $E$
(with $u_i=I_{E_i}$), this implies that for any biorthogonal system
$(x_1,\ldots,
x_n)$ $(x^*_1,\ldots, x^*_n)$ in $E$ we have
$$n \le \left\|\sum^n_1 \lambda(g_i) \otimes x_i\right\|_{C_\lambda
\otimes_{\rm min} E} \left\|\sum \lambda (g_i)\otimes
x^*_i\right\|_{C_\lambda\otimes_{\rm min} E^*}.\leqno (2.5)$$
Now let $E^\lambda_n = \hbox{span}(\lambda(g_1),\ldots, \lambda(g_n))
\subset C_\lambda$.

\n Let $\lambda_*(g_i)$ be the biorthogonal functionals
in $(E^\lambda_n)^*$. Then, if $n>1$, we have by [AO]
$$\left\|\sum^n_1\lambda(g_i) \otimes \lambda(g_i)\right\| =
2\sqrt{n-1},$$
and since $t =  \sum\limits^n_{i=1} \lambda(g_i) \otimes
\lambda_*(g_i)$
represents the inclusion map $j\colon \ E^\lambda_n\to C_\lambda$, we
have
$$\|t\|_{C_\lambda\otimes_{\rm min} (E^\lambda_n)^*} =
\|j\|_{cb(E^\lambda_n, C_\lambda)} =1.$$
Hence taking $x_i = \lambda(g_i), x^*_i = \lambda_*(g_i)$ in (2.5) we
obtain
$$n\le 2\sqrt{n-1}$$
or equivalently $n\le 2$.\qed

\n {\bf Remark 2.4.} By a simple modification of the preceding proof,
one can prove that the subset $HOS_n\subset OS_n$ formed of all the
$n$-dimensional operator spaces which are isometric to $\ell_2^n$
is non-separable if $n>2$. Our original argument here gave
only $n>4$, the improvement is due to Timur Oikhberg.

\n Here is briefly the argument: By the proof
of Theorem 2.3, if $HOS_n$ is
separable, then any $E$ in $HOS_n$ must satisfy (2.5).
 Consider then the operator space
$\min(\ell_2^n)$ obtained by embedding $\ell_2^n$
isometrically
into a commutative $C^*-$algebra. Let $(e_i)$ be the
 basis of
$\ell_2^n$. Assume $n>1$.
Let $\gamma_n=2(1-n^{-1})^{1/2}$. Consider the subspace
$$E\subset E^\lambda_n\oplus \min(\ell_2^n)$$
spanned by the vectors $x_i=\lambda(g_i)\oplus \gamma_n e_i$.
Then by [AO, p. 1038] we have $\|\sum \alpha_i x_i\|=
\gamma_n(\sum|\alpha_i|^2)^{1/2},\break \forall (\alpha_i)\in \comp^n$,
hence $E\in HOS_n$. Furthermore, we have
$$\left\|\sum^n_1 \lambda(g_i) \otimes x_i\right\|_{C_\lambda
\otimes_{\rm min} E} =\max\{2\sqrt{n-1},\gamma_n^2\}\quad{\rm and}\quad
\left\|\sum \lambda (g_i)\otimes
x^*_i\right\|_{C_\lambda\otimes_{\rm min} E^*}\le 1,\leqno(2.6)$$
so that we conclude again from (2.5) that $n\le 2$.

\n{\bf Remark 2.5.} (i) An operator space $E$ is called homogeneous
if for any $u\colon E\to E$ we have $\|u\|=\|u\|_{cb}$. This notion
seems particularly interesting in the Hilbertian case (see [P3]).
Consider an arbitrary $n$-dimensional operator space $E$
given with a basis   $(e_1,...,e_n)$. We can define its "homogeneous
hull" $\hat E$ as follows. Let $U(n)$ be the unitary group. We view
the coordinates
$(u_{ij})$ of a unitary matrix $u$ as a continuous function
on $U(n)$, so that $u_{ij}\in C(U(n))$. In the space
$C(U(n))\otimes_{\min} E$ we consider the elements
$$\hat e_i=\sum_j u_{ij} \otimes e_j \in C(U(n))\otimes_{\min} E.$$
Let $\hat E$ be the operator space spanned by
$(\hat e_1,...,\hat e_n)$. Observe that $U(n)$
acts isometrically (and actually completely isometrically)
on $\hat E$, therefore it is easy to check that $\hat E$ is Hilbertian
and homogeneous. Moreover,
we have $d_{SK}(\hat E)\le d_{SK}(E)$.
Now let $F$ be another operator space
and let $u\colon E\to F$ be an isomorphism. Let
$f_i=u(e_i)$ and let $\hat F$ be the operator  space associated to F
and this basis. Then it is easy to check that
$$d_{cb}(\hat E,\hat F)\le \|u\|_{cb} \|u^{-1}\|_{cb}
 \quad {\rm and}
 \quad d_k(\hat E,\hat F)\le \|u\|_{k} \|u^{-1}\|_{k}\ \ \forall k.$$

\n (ii) Let us denote by $HH_n$ the subset of $OS_n$ formed
of all the Hilbertian homogeneous
spaces.
Using the first part of this remark, it is easy to check that any
space $E$ in $HH_n$ is the weak limit
of a net $(E_i)$ in $HH_n$ such that $d_{SK}(E_i)=1$ for all $i$.
Then a simple modification of the proof of Theorem 2.3 shows that
$HH_n$ is a non-separable subset of $OS_n$
 for the (strong) distance $\delta_{cb}$, if $n>2$.
Indeed, we can replace the space $E$ and its basis $(x_i)$ in Remark
2.4
by $\hat E$ and $\hat x_i$. Using an inequality
due to Haagerup [H3, Lemma 2.4], one can check that (2.6)
 remains valid. This gives us a space in $HH_n$
which satisfies (2.5) only if $n\le 2$. (Alternately, one could
replace $E$ by the linear span of a circular system
in the sense
of Voiculescu [VDN], but this seems to yield
 non-separability only for  $n>4$.)

\proclaim Proposition 2.6.
 Let $A$ be any separable $C^*$-algebra
 (or any separable operator
space) and let $n\ge 1$.
\item{\rm (i)} The subset $S_n(A)\subset OS_n$ formed of
all the $n$-dimensional subspaces of $A$ is
separable in $(OS_n,\delta_{cb})$.
\item{\rm (ii)} For any $n>2$,
 there is an operator space
$E_0\in OS_n$ and $\varepsilon_0>0$
such that for any $n$-dimensional subspace $E\subset A$ we have
$d_{cb}(E,E_0)\ge 1+\varepsilon_0$.

\pf The first part is proved by a standard perturbation
argument.
We merely sketch it: Let $D\subset A$ be a dense
countable subset. For any $n$-tuple $x=(x_1,...,x_n)$
of linearly independent elements of $A$, let $E_x\subset
A$ be their linear span. Let $x_1^*,...,x_n^*$ be
functionals in $A^*$ which are biorthogonal to $x_i$. Fix
$\vp>0$. Pick $y_1,...,y_n$ in $D$ such that
$\sum \|x_i-y_i\|<\vp$. Consider then the operator
$u\colon \ A\to A$ associated to $\sum x_i^*\otimes
(y_i-x_i)$. Clearly $\|u\|_{cb}\le f(\vp)$ with $f(\vp)\to
0$ when $\vp\to 0$. Moreover,
$(I+u)(x_i)=y_i$, and (say, when
$f(\vp)<1$)  $x_i=(I+u)^{-1}(y_i)$. This immediately yields
$d_{cb}(E_x,E_y)\le 1+g(\vp)$ with $g(\vp)\to 0$ when
$\vp\to 0$, whence (i).

\n The second part follows from the first one, since it
merely expresses the fact that (by Theorem
2.3) $S_n(A)$ is not dense in
$OS_n$ for $n>2$.
 \qed
We now give a more precise version of Theorem~2.3 based on the
following
well known variant of Baire's Theorem.

\proclaim Lemma 2.7. Let $S,T$ be metric spaces and let $f\colon \ S\to
T$
be a mapping such that for every closed ball $B\subset T$, $f^{-1}(B)$
is
closed in $S$. Fix a number  $\varepsilon>0$. Let $C_\varepsilon(f)$ be
the
set of all points of $\vp$-continuity of $f$, i.e.\ all points $s$ in
$S$
such that whenever $s_i\to s$ we have $\limsup d(f(s_i), f(s)) \le
\vp$.
Now assume that $T$ is $\vp$-separable, i.e.\ there exists a sequence
$\{B_n\}$ of closed balls of radius $\vp$ in $T$ such that $T = \cup
B_n$.
Then, if $S$ is a Baire space, the set $C_{2\vp}(f)$ of points of
$(2\vp)$-continuity of $f$ is dense in $S$.

\n {\bf Proof.} Let $S' = \bigcup\limits_n [f^{-1}(B_n) \backslash
\mathop{f^{-1}(B_n)}\limits^{\buildrel \circ\over \frown}] \subset S$.
Since $S'$ is a countable
union of closed sets with empty interior in a Baire space,
its complement $S\backslash S'$ is dense in $S$. But it is
easy to check (recall $S = \cup f^{-1}(B_n)$), that
$S\backslash S'\subset C_{2\vp}(f)$.\qed

\proclaim Theorem 2.8. For each $n\ge 1$, let $\vp_n$ be the infimum of
all
numbers $\vp>0$ such that $(OS_n, \delta_{cb})$ is
$\vp$-separable. Let $\delta_n = \exp(\vp_n)$. Then we
have for all $n\ge 3$ $$\left({n\over
2\sqrt{n-1}}\right)^{1/4} \le\delta_n.\leqno (2.7)$$ and
there is a constant $c>0$ such that for all $n\ge 3$
$$cn^{1/8} < \delta_n \le n^{1/2}.\leqno (2.8)$$

\n {\bf Proof.} We apply Lemma~2.7 to the same map $f$ as in the proof
of
Theorem~2.3.

\n Consider $\delta>\delta_n$ and let $\vp = \log(\delta)$
so that $(OS_n, \delta_{cb})$ is $\vp$-separable.
By Lemma 2.7, $C_{2\vp}(f)$ is $\delta_w$-dense. Then, by a
simple modification of the proof of Proposition~2.2,
 any $E$ is $OS_n$
is the weak limit of a net $E_i$ such that $\limsup
d_{SK}(E_i) \le e^{2\vp}= \delta^2$ and $\limsup
d_{SK}(E^*_i) \le e^{2\vp} = \delta^2$. Then (recalling
Corollary~1.10) we conclude, as in the proof of
Theorem~2.3, that $$n\le \delta^4 2\sqrt{n-1}$$ whence
(2.7).

\n It was proved in [P3, Theorem 9.6] that $d_{SK}(E)\le
n^{1/2}$ for all $E$ in $OS_n$, therefore $\delta_n \le
n^{1/2}$ and (2.8) follows.\qed

\proclaim Corollary 2.9. There is a constant $c>0$ such that for each
$n\ge
3$ there is an uncountable  collection $(E_i)_{i\in I}$ is $OS_n$
satisfying
$$\forall\ i\ne j\qquad d_{cb}(E_i,E_j) > cn^{1/8}$$

It would be interesting to find the exact
asymptotic behaviour of $\delta_n$.

\n {\bf Remark 2.10.} Two or three months after this paper
had been circulated as a preprint, Simon Wassermann
mentioned to us that he conjectured that the linear spans
of the $n$-tuples of operators considered by Voiculescu  in
[V] should yield a more explicit non-separable
family of finite dimensional operator spaces.
Some form of this conjecture is indeed correct. This shows that
groups with Kazhdan's property T (see [DHV]) can be used to prove the
non-separablility
of $OS_n$. Here are the details.
In [V], to each subset $\Omega$ of the integers, Voiculescu associates
an $n$-tuple  $\tau_\Omega=(T^\Omega_1,...,T^\Omega_n)$
 of operators in $B(H)$ (with say $H=\ell_2$), in the following way.
Let $G$ be any discrete group with Kazhdan's property T
admitting a countable collection of pairwise disjoint
finite dimensional representations $(\pi_k)$. (For instance
we can take $G=SL_3(\ent)$.) Let $(g_1,...,g_{n-1})$ be a finite set
of generators and let $g_n$ be equal to the unit element. Then
for any $\Omega\subset \nat$, we define
$$T^\Omega_j=\oplus_{k\in \Omega} \pi_k(g_j),\quad j=1,2,...,n.$$
Let $M_\Omega$ be the finite dimensional operator space
spanned by $\tau_\Omega$.
We claim these form a non-separable collection of
operator spaces.
This does not seem to follow from Voiculescu's stated results
but it does follow easily from the proof of his proposition 3.
Indeed  he shows a stronger
result than he states, as follows:

\n First, in case $M_\Omega$ is not $n$-dimensional
we consider an $n$-dimensional
operator space $E_\Omega$ containing $\tau_\Omega$.

\n Following [V], we use property T through the
following: There is a fixed number
$\vp>0$ such that if $\Omega$ and $\Omega'$ are subsets of
$\nat$ with  $\Omega\not\subset\Omega'$, then
there are unitary operators $u_1,...,u_n$
satisfying
$$  n=\|\sum_1^n T^\Omega_j\otimes u_j\|_{\min}\quad{\rm and}\quad
\|\sum_1^n T^{\Omega'}_j\otimes u_j\|_{\min}< n-\vp .\leqno(2.9)$$
(More precisely, if we pick $k\in \Omega-\Omega'$
then we can take $u_j=\overline{\pi_k(g_j)}$.)
\def\c{{\cal C}}

\n Now assume that the metric space $(OS_n, \delta_{cb})$
 is separable. Let
$(E_m)$ be a dense  sequence in $OS_n$.
Fix a number $\eta>0$. Then there is an integer $m$ and a continuous
collection $\c$ of subsets of $\nat$
such that for each $\Omega$ in $\c$ there is a map
$v_\Omega\colon E_\Omega\to E_m$ such that
$\|v_\Omega\|_{cb}<1+\eta$ and $\|v^{-1}_\Omega\|_{cb}=1$.
Now consider the continuous family
$(v_\Omega(T^\Omega_j))_{j\le n}$ of $n$-tuples of elements of $E_m$.
Since $(E_m)^n$ (the space of $n$-tuples of elements of $E_m$)
is norm-separable, there must exist
a continuous subcollection $\c_1\subset\c$ such that
for all $\Omega,\Omega'$ in $\c_1$ we have
$$ \sum_1^n\|v_\Omega(T^\Omega_j)-v_{\Omega'}(T^{\Omega'}_j)\|<\eta.
\leqno(2.10)$$
A fortiori, $\c_1$ has cardinality $>1$, hence we can find
  $\Omega,\Omega'$ in $\c_1$
 satisfying $\Omega\not\subset{{\Omega'}}$.
By (2.9) this implies
$$\|\sum v_{{\Omega'}}(T^{{\Omega'}}_j)\otimes u_j\|_{\min}
\le \|v_{\Omega'}\|_{cb} \|\sum T^{\Omega'}_j\otimes u_j\|_{\min}
\le (1+\eta)(n-\vp),$$
and $$n=(\|v^{-1}_\Omega\|_{cb})^{-1}\|\sum_1^n T^\Omega_j\otimes
u_j\|_{\min}
\le \|\sum v_\Omega (T^\Omega_j)\otimes u_j\|_{\min}.$$
This gives by (2.10)
$$n\le \| \sum v_{\Omega'}(T^{\Omega'}_j)\otimes u_j\|_{\min}
+\sum_1^n\|v_{\Omega'}(T^{\Omega'}_j)-v_\Omega(T^\Omega_j)\|
\le (1+\eta)(n-\vp)+\eta.$$
When $\eta>0$ is small enough this is impossible
(since $\vp>0$ remains fixed). This contradiction
completes the proof.
\qed

\n Actually, Simon  Wassermann conjectured that
the whole
family of spaces $(M_\Omega)$ is uniformly
 $d_{cb}$-separated, i. e. that for some
$\eta>0$ we have $d_{cb}(M_\Omega,M_{\Omega'})>1+\eta$
whenever $\Omega\not=\Omega'$. As far as we know this
is still open.

\n It is known ([Tr]) that $G=SL_3(\ent)$ admits two generators
so that (recall that $g_n$ is the unit) we  obtain by
 this reasoning a
continuous
collection of $3$-dimensional operator spaces $(E_t)$
such that for some $\vp>0$ we have
$d_{cb}(E_t,E_s)>1+\vp$ for all $t\ne s$. By a simple
modification, we can
make sure
that the spaces we obtain  are spanned by
three unitaries. The same cannot be achieved
with spans of two   unitaries. Indeed, the span of
two unitaries $(U_1,U_2)$ is completely isometric
to the span of $(I,U)$ with $U=U_1^*U_2$, which itself
embeds completely isometrically into a commutative
(hence nuclear) $C^*$-algebra.

\n However, this approach does not seem to yield the refinements
 in Remarks  2.4 and 2.5 or in Corollary 2.9.

\vfill\eject

\n {\bf \S 3. Applications to $B(H)\otimes B(H)$ and maximal operator
spaces.}

By Kirchberg's results in [Ki2], Theorem 2.3 implies

\proclaim Corollary 3.1. If $\dim H = \infty$, there is more than one
$C^*$-norm on $B(H)\otimes B(H)$. In other words we have
$$B(H)\otimes_{\rm min} B(H)\ne B(H)\otimes_{\rm max} B(H).$$

\pf This follows from the equivalence of the conjectures (A7) and (A2)
in [Ki2, p.
483]. For the convenience of the reader,
we include a direct argument, as
follows. For any discrete group $G$
we denote by $C^*(G)$ the full $C^*$-algebra of $G$.
 By Proposition 2.6 (ii) (applied with
$A=C^*(F_\infty)$), there is  an operator space
 $E_0$  such that for some $\varepsilon_0>0$ we
have
$$d_{cb}(E,E_0) \ge 1+\varepsilon_0$$
for all $n$-dimensional subspaces $E\subset C^*(F_\infty)$.
 Now let $F_I$ be a free group associated to a set of generators
 $\{g_i\mid
 i\in I\}$ where $I$ is any set with infinite cardinality. Observe that
 for
 any finite dimensional (or merely separable) subspace $E\subset
 C^*(F)$
 there is a countable infinite subset $J\subset I$ such that $E\subset
 C^*(F_J)$. (Indeed only countably many ``letters'' are being used.)
 Hence
 we also have
 $$d(E,E_0) \ge 1+\vp_0$$
 for all $n$-dimensional subspaces $E\subset C^*(F_I)$.
Let $\pi\colon \ C^*(F)\to B(H)$ be a $C^*$-algebra representation of
the
full $C^*$-algebra of a big enough free group $F$ onto $B(H)$.
 We have
 $$B(H)\approx C^*(F)/\hbox{Ker } \pi.$$

\n Let ${\cal I} = \hbox{Ker } \pi$. Then (by Lemma 0.8)
the quotient norm of the space $$Q =
[C^*(F)\otimes_{\rm min} B(H)]/[{\cal I}\otimes_{\rm min} B(H)]$$
induces on
$B(H)\otimes B(H)$ a $C^*$-norm. Assume that there is only one such
norm.
Then we have isometrically
$$Q = B(H)\otimes_{\rm min} B(H).\leqno(3.1)$$
Now we may clearly assume that
$E^*_0$ is a subspace of $B(H)$, i.e.\ we have
$E^*_0\subset B(H)$ completely isometrically, so
 that the completely isometric inclusion
$j\colon\ E_0\to B(H)$ can be viewed as an
 element $j_0$ in $B(H)
\otimes_{\rm min} E^*_0\subset B(H)\otimes_{\rm min} B(H)$.
 Note that
$\|j_0\|_{\rm min} = \|j\|_{cb}=1$.
 By Lemma~0.7 and by (3.1),
for any $\varepsilon>0$
there is a lifting $\tilde j_0$ in $C^*(F)\otimes_{\rm min}E^*_0 $ with
$\|\tilde j_0\|_{\rm min} < 1+\varepsilon$ and $(\pi \otimes I_{E^*_0})
(\tilde j_0)
= j_0$. Now, let $\tilde j\colon \ E_0\to C^*(F)$ be the associated
linear
operator. We have $\pi\tilde j = j$ and
$\|\tilde j_0\|_{\rm min} = \|\tilde j\|_{cb}$, hence $d_{cb}(E_0,
\tilde j(E_0)) \le
\|\pi\|_{cb} \|\tilde j\|_{cb}<1+\varepsilon$.
When $\varepsilon<\varepsilon_0$ this is impossible.\qed

\n{\bf Remark.} By [Ki2, section 8], Corollary 3.1 has the following
consequences

(i) There is a separable unital $C^*$-algebra with the WEP in the sense
of
Lance [La] for which Ext$(A)$ is not a group.

(ii) There are separable unital $C^*$-algebras $A,B$ with WEP such that
$$A\otimes_{\min}B\not=A\otimes_{\max}B.$$

(iii) The WEP does not imply the local lifting property (in short LLP)
in the
sense of [Ki2].

(iv) There is a separable unital $C^*$-algebra $A$ with WEP which is
not
approximately injective in the sense of [EH].

(v) The identity $A^{op}\otimes_{\min}A=A^{op}\otimes_{\max}A$ does not
imply the
approximate injectivity of $A$.

(vi) There is a unital separable $C^*$-algebra $B$ with WEP which is
not a
quotient $C^*$-algebra of an approximately injective $C^*$-algebra.

By [Ki2, p. 484], Corollary 3.1 also implies a negative answer to
Kirchberg's
conjecture (C1) in [Ki2]. Thus we have :

$C^*(F_\infty)$ and $C^*(SL_2(\ent))$ are not approximately injective
in $B(H)$.

Finally, by [Ki2, p. 487], Corollary 3.1 implies a negative answer to
the
conjecture (P2) in [Ki2], hence we have : Let $K=K(H),\ B=B(H)$.
Consider the
canonical morphism
$$\Phi\colon B\otimes_{\min}B\mapsto(B/K)\otimes_{\min}(B/K).$$
Then the kernel of $\Phi$ is strictly larger than the set
$F(K,B,B\otimes_{\min}B)+F(B,K,B\otimes_{\min}B)$ where $F(.,.,.)$
denotes the
Fubini product.

In a different direction, we can give more examples of non-exact
operator
spaces, completing those of [P1]. Following Blecher and
Paulsen [BP], given a Banach
space $E$, we denote by $\min(E)$ the operator space obtained by
embedding
$E$ into a commutative $C^*$-algebra, or equivalently by embedding $E$
into
the space $C(K)$ of all continuous functions on $K$ with $K = (B_{E^*},
\sigma(E^*,E))$. Similarly, let $I$ be a suitable set and let $g\colon
\
\ell_1(I)\to E$ be a metric surjection (i.e.\ $q^*$ is an isometry). We
view the space $\ell_1(I) = c_0(I)^*$ as an operator space with the
dual
operator space structure. Then (cf.\ [BP, Pa2]) we denote by $\max(E)$
the
operator space obtained by equipping $E$ with the operator space
structure
(in short o.s.s.) of the quotient space $\ell_1(I)/\hbox{Ker}(q)$.
Equivalently, we have a complete isometry
$$\max(E)\longrightarrow (\min(E^*))^*,$$
and $\max(E)$ is characterized by the isometric identity
$$cb(\max(E), M_n) = B(E,M_n).\leqno (3.2)$$
More generally, for any operator space $F$ we have isometrically
$$cb(\max(E), F) = B(E,F).\leqno  (3.2)'$$
We refer to [BP, Pa2] for more information. It will be convenient to
introduce
the following characteristic for an operator space $E$.
$$d_{QSK}(E) = \inf\{d_{cb}(E,F)\}$$
where the infimum runs over all operator spaces $F$ of the form $F =
E_1/E_2$ where $E_2\subset E_1\subset K$ and $\dim E = \dim F$.

\proclaim Theorem 3.2. Let $E$ be any $n$-dimensional Banach space.
Then
$$d_{QSK}(\max(E) )\ge {\sqrt n\over 4}.$$

\pf By a well known result in Banach space theory
 (cf.\ e.g.\ [P2, Theorem 1.11 p.15]) we have
$\pi_2(I_{E^*}) = \sqrt n$.
Let $u\colon\ X\to Y$ be a 2-absolutely summing operator between Banach
 spaces and let  $J\colon \ Y\to Y_1$ be an isometric embedding. Then
 it is easy to see that
 $$\pi_2(u) = \pi_2(Ju).\leqno (3.3)$$
 Consider an isomorphism $v\colon \ E_1/E_2\to \max(E)$ with
 $E_2\subset
 E_1\subset K$. Observe
 $(E_1/E_2)^* = E^\bot_2\subset E^*_1$. Let $J\colon \ E^\bot_2 \to
 E^*_1$
 be the (isometric) canonical inclusion. We will now apply
 Corollary~1.7 to
 the map
 $$Jv^*\colon \ \min(E^*) {\buildrel v^*\over \longrightarrow} E^\bot_2
 {\buildrel J\over \longrightarrow} E^*_1.$$
 Note that $d_{SK}(\min(E^*)) = d_{SK}(E_1) = 1$. Hence Corollary~1.7
 yields
 by (3.3)
 $$\pi_2(v^*) = \pi_2(Jv^*) \le 4\|Jv^*\|_{cb} \le 4\|v^*\|_{cb} =
 4\|v\|_{cb}$$
 but this implies by (0.4)
 $$n^{1/2} = \pi_2(I_{E^*}) \le \|v^{-1*}\| \pi_2(v^*) \le 4\|v^{-1}\|
 \
 \|v\|_{cb}$$
 hence $n^{1/2} \le 4d_{QSK}(\max(E))$.\qed

\n {\bf Remark.} In particular the preceding result answers a question
raised by Vern Paulsen (private communication): \ the space
$\max(\ell^n_2)$ is quite different (when $n$ is large enough) from the
linear span of $n$ Clifford matrices, i.e.\ matrices $(u_i)$ in
$M_{2^n}$
satisfying the relations
$$u_i=u^*_i\qquad u_iu^*_j + u^*_ju_i = 2\delta_{ij}I.$$
Indeed, if we denote $Cl_n = \hbox{span}(u_1,\ldots, u_n)$ then by
Theorem~3.2 we have
$$d_{cb}(\max(\ell^n_2), Cl_n) \ge d_{QSK}(\max(\ell^n_2)) \ge \sqrt
n/4.$$

\vfill\eject
\magnification\magstep1
\baselineskip = 18pt
\def\n{\noindent}
\def\n{\noindent}
\def\vp{\varepsilon}
\def\cf{{\it cf.\/}\ }
\def\ent{{{\rm Z}\mkern-5.5mu{\rm Z}}}
 \overfullrule = 0pt
\def\C{C^*(F_\infty)}
 \def\qed{{\hfill{\vrule height7pt width7pt
depth0pt}\par\bigskip}}
\def\implies{\Rightarrow}

\let\eps =\varepsilon
\def \comp{ \;{}^{ {}_\vert }\!\!\!{\rm C}   }
\def\N#1{\left\Vert#1\right\Vert}
\def \nat{ { {\rm I}\!{\rm N}} }

\overfullrule = 0pt
\def\pf{\medskip{\noindent{\bf Proof. }}}

\n {\bf \S 4. A new tensor product for $C^*$-algebras or
operator spaces.}

 Let $E_1\subset B(H_1),E_2\subset B(H_2)$ be arbitrary operator
 spaces.

\n Let us denote by $\N{~}_M$ the norm induced on the algebraic tensor
product $E_1\otimes E_2$ by
$B(H_1)\otimes_{\max}B(H_2)$. Recalling (0.3), and using
the injectivity of $B(H_1)$ and $B(H_2)$ as well as the
decomposition property of $c.b.$ maps into $B(H)$, it is easy
to check  (see
Lemma 4.1 below) that this norm is independent of the
choice of the completely isometric embeddings $E_1\subset
B(H_1),E_2\subset B(H_2)$. In other words, the norm
$\N{~}_M$ on $E_1\otimes E_2$ depends only on the operator
space structures of $E_1$ and $E_2$.

\n  We will denote by
$E_1\otimes_M E_2$ the completion of $E_1\otimes E_2$ under this norm.
We equip the space
$E_1\otimes_M E_2$ with the natural operator space structure induced by
the $C^*$-algebra
$B(H_1)\otimes_{\max} B(H_2)$ via the isometric embedding
 $E_1\otimes_M E_2\subset B(H_1)\otimes_{\max}B(H_2)$.

\n Clearly, if $A,\  B$ are $C^*$-algebras, then
$\|\ \|_M$ is a $C^*$-norm on $A\otimes  B$ and
$A\otimes_M B$ also is a $C^*$-algebra.

\proclaim Lemma 4.1.  Let $F_1,F_2$ be two operator spaces.
Consider $c.b.$ maps $u_1:E_1\to F_1$ and $u_2:E_2\to
F_2$.
 Then $u_1\otimes u_2$ defines a $c.b.$ map from
$E_1\otimes_M E_2$ to $F_1\otimes_M F_2$ with
$$\N{u_1\otimes u_2}_{{cb(E_1\otimes_M E_2,F_1\otimes_M F_2)}}\leq
\N{u_1}_{cb}\N{u_2}_{cb}.\leqno(4.1)$$

\pf
Indeed, note that if $F_1\subset B({\cal H}_1),\ F_2\subset B({\cal
H}_2)$ then by the extension
property of $c.b.$ maps (\cf [Pa1] p. 100)  $u_1,u_2$ admit extensions
$\tilde u_1:B(H_1)\to
B({\cal H}_1)$ and $\tilde u_2:B(H_2)\to B({\cal H}_2)$ with $\N{\tilde
u_1}_{cb}=\N{u_1}_{cb}$
and $\N{\tilde u_2}_{cb}=\N{u_2}_{cb}$. Hence it suffices to check this
in the case when each of
$E_1,E_2,F_1,F_2$ is $B(H)$ for some $H$.
Then the idea is to use
the decomposition property of $c.b.$ maps on $B(H)$ as
linear combinations of completely positive maps to reduce
checking (4.1) to the case of  completely positive maps.
In the completely positive case,
the  relevant point here is of course (0.3). This idea
leads to a simple proof of  (4.1) with some additional
numerical factor. However, this factor can be removed at
the cost of a slightly more technical argument based on
[H3]. For lack of a suitable reference,  we now briefly
outline this (straightforward) argument to check (4.1).

\n Let
$A,B$ be $C^*$-algebras, we will denote by $CP(A,B)$ the
set of all completely positive maps $u\colon \ A\to B$,
and by ${D}(A,B)$ the set of all decomposable maps $u\colon
\ A\to B$, {\it i.e. } maps which can be written as $u =
u_1-u_2 + i(u_3-u_4)$ with $u_1,\ldots, u_4 \in CP(A,B)$.

\n In [H3], Haagerup defines the norm $\|u\|_{dec}$ on
${D}(A,B)$ as follows. Consider all possible mappings
$S_1,S_2$ in $CP(A,B)$ such that the map $v\colon \ A\to
M_2(B)$ defined by $$v(x) =
\left(\matrix{S_1(x)&u(x^*)^*\cr\cr
u(x)&S_2(x)\cr}\right)$$ is completely positive. Then we
set $\|u\|_{dec} = \inf\{\max\{\|S_1\|, \|S_2\|\}\}$ where
the infimum runs over all possible such mappings.

\n In [H3, Proposition~1.3 and Theorem~1.6] the following
results appear: $$\leqalignno{\forall\ u\in {D}(A,B)\qquad
\|u\|_{cb} &\le \|u\|_{dec},&(4.2)\cr
\forall\ u\in CP(A,B)\qquad \|u\| &= \|u\|_{cb} =
\|u\|_{dec},&(4.3)}$$
 (see also [Pa1, p.~28] for the first equality)
$$\hbox{if $C$ is any $C^*$-algebra,
if $u\in D(A,B)$ and $v\in D(B,C)$, then $vu\in D(A,C)$
and}\leqno (4.4)$$ $$ \|vu\|_{dec}
\le \|v\|_{dec} \|u\|_{dec},$$ $$\forall\ u \in cb(B(H),
B({\cal H}))\qquad \|u\|_{cb} = \|u\|_{dec}. \leqno (4.5)$$

\n Now let $A_1,A_2,B_1,B_2$ be arbitrary $C^*$-algebras,
and let $u_1 \in {D}(A_1,B_1)$, $u_2\in {D}(A_2,B_2)$. We
claim that Haagerup's results imply that $u_1\otimes u_2$
extends to a decomposable map from $A_1\otimes_{\rm max}
A_2$ into $B_1\otimes_{\rm max}B_2$ (still denoted by
$u_1\otimes u_2$) satisfying $$\|u_1\otimes u_2\|_{dec}
\le \|u_1\|_{dec} \|u_2\|_{dec}.\leqno (4.6)$$ By (4.2) we
have a fortiori $$\|u_1\otimes u_2\|_{cb(A_1\otimes_{\rm
max} A_2, B_1\otimes_{\rm max}B_2)} \le \|u_1\|_{dec}
\|u_2\|_{dec}.\leqno (4.7)$$ To verify (4.6) (hence also
(4.7)) we may assume (using (4.4) and $u_1\otimes u_2 =
(u_1\otimes I) (I\otimes u_2)$) that $A_2=B_2$ and $u_2$
is the identity on $A_2$. Consider then $v_1\colon \
A_1\to M_2(B_1)$ of the form $$v_1(x) =
\left(\matrix{S_1(x)&u_1(x^*)^*\cr\cr
u_1(x)&S_2(x)\cr}\right)$$ with $v_1$ and $S_1,S_2$ all
completely positive.

\n By (0.3),  the associated map $v_1\otimes
I_{A_2}$
as well as
$S_1\otimes I_{A_2}$ and $S_2\otimes I_{A_2}$
are (bounded and) completely
positive
from $  A_1\otimes_{\rm max} A_2$ to
$B_1\otimes_{\rm max} A_2$. Therefore (by the definition
of $\|\ \|_{dec}$) $$\leqalignno{\|u_1\otimes
I_{A_2}\|_{dec} &\le \max\{\|S_1\otimes I_{A_2}\|,
\|S_2\otimes I_{A_2}\|\}\cr
 &\le \max\{\|S_1\|, \|S_2\|\}.&\hbox{hence by
(0.3)}}$$ It follows that $$\|u_1\otimes I_{A_2}\|_{dec}
\le \|u_1\|_{dec},$$ and this is enough to verify (4.6)
(and a fortiori (4.7)).

\n The proof of (4.1) is now easy:\ let $u_1,\widetilde
u_1$ and $u_2,\widetilde u_2$ be as explained above, by
(4.7) $$\leqalignno{\|u_1\otimes u_2\|_{cb(E_1\otimes_M
E_2, F_1 \otimes_MF_2)} &\le \|\widetilde u_1\otimes
\widetilde u_2\|_{cb(B(H_1)\otimes_{\rm max} B(H_2),
B({\cal H}_1) \otimes_{\rm max} B({\cal H}_2))}\cr & \le
\|\widetilde u_1\|_{dec} \|\widetilde u_2\|_{dec}\cr &\le
\|\widetilde u_1\|_{cb} \|\widetilde u_2\|_{cb} =
\|u_1\|_{cb} \|u_2\|_{cb}.&\hbox{hence by (4.5)}}$$ This
completes the proof of (4.1). \qed

\n We will   use several times the following
 obvious consequence of (4.1) and the definition of
$E_1\otimes_M E_2$:

\n (4.8) \ If $u_1$ and $u_2$ (as above) are complete isometries, then
$u_1\otimes u_2\colon E_1\otimes_M E_2\to F_1\otimes_M F_2$ also is a
complete isometry.

By Corollary 3.1, we know that there are operator spaces $E,F$ such
that
$E\otimes_{\min}F\not=E\otimes_M F$. It is natural to try to understand
the meaning of
this new tensor norm $\|\ \|_M$ and to characterize the operator
spaces
$E,F$ for which the equality holds. For that purpose the following
result due to Kirchberg [Ki3] will be crucial:

\n For any free group $F_I$ and any $H$, we have an isometric identity
$$    C^*(F_I)\otimes_{\min} B(H) =  C^*(F_I)\otimes_{\max} B(H)
.\leqno(4.9)$$
Using this, we have

\proclaim Proposition 4.2. Let $E,F$ be operator spaces, let $u\in
E\otimes F$
and let $U\colon F^*\to E$ be the associated finite rank linear
operator.
Consider a finite dimensional subspace $S\subset \C$ and
a factorization of $U$ of the form  $U=ba$
with bounded linear maps $a\colon F^*\to S$ and $b\colon S\to E$, where
$a\colon F^*\to S$ is weak-$*$ continuous. Then
$$\|u\|_M=\inf \{ \|a\|_{cb} \|b\|_{cb}\}$$
where the infimum runs over all such factorizations
of $U$.

\pf Assume $E\subset B(H)$
and $F\subset B(K)$ with $H,K$ Hilbert.
It clearly suffices to prove this in the case when $E$ and
$F$ are both finite dimensional. Assume $U$ factorized as above
with  $ \|a\|_{cb} \|b\|_{cb}< 1$.
Then by Kirchberg's theorem (4.9) the $\min$ and $\max$ norms are equal
on $\C\otimes B(K)$, hence, by (4.8), we have isometrically
$S\otimes_{\min} F=S\otimes_{M} F$, so that if $\hat a$ is the element
of $S\otimes_{\min} F$ associated to $a$, we have
$\|\hat a\|_M=\|a\|_{cb}$ and
$u= (b\otimes I_F)(\hat a)$. Therefore,
by (4.1) we have $\|u\|_M\le \|b\|_{cb} \|\hat a\|_M
\le\|a\|_{cb} \|b\|_{cb}<1$.

\n The proof of the converse is essentially the same as for
Corollary 3.1 above.  We skip the details.\qed

 Let
$A$ be any $C^*$-algebra. For any finite dimensional operator space
$E$, let
$$d_{SA}(E)=\inf\{d_{cb}(E,F)|\ F\subset A\}.$$
Then the preceding result immediately implies
\proclaim Corollary 4.3. Let $E$ be a finite dimensional operator
space.
Let $i_E\in E\otimes E^*$ be the tensor associated to the identity on
$E$.
Then
$$\|i_E\|_M= d_{S\C}(E).$$
In particular we have
$$d_{S\C}(E)= d_{S\C}(E^*).\leqno(4.10)$$

\pf The first part is clear by Proposition 4.2. To check (4.10),
observe more generally that for any operator spaces
$E,F$, the "flip isomorphism" ($x\otimes y \to y\otimes x$)
is a complete isometry between
 the spaces $E\otimes_M F$ and $F\otimes_M E$. Hence (4.10)
 follows
by symmetry. \qed

\n{\bf Remark.}
The preceding argument shows the following : If
$A\otimes_{\min}B(H)=A\otimes_{\max}B(H)$ then for any finite
dimensional
operator space $E$ we have
$$d_{{SC^*}(F_\infty)}(E)\leq d_{SA}(E^*).\leqno(4.11) $$
In particular,
this holds (by definition of nuclearity)
if $A$ is nuclear. (Hence (4.11) still holds if $A$
is exact, since $ d_{SK}(E^*)\le  d_{SA}(E^*)$
in that case).
The proof of Corollary 3.1 shows that (4.11) is false
in general for $A=B(H)$. (Observe that
  $d_{SA}(E^*)=1$    for any $E$ if $A=B(H)$ and $\dim(H)=\infty$.)

\n {\bf Remark.} We can now give a quantitative version of
Corollary~3.1. For any $n$ let
$$\lambda(n) =
\sup\left\{{\|u\|_{\rm max}\over \|u\|_{\rm
min}}\right\}\leqno (4.12)$$
where the supremum runs over all $u$ in $B(H)\otimes B(H)$ with
$\hbox{rank}(u)\le n$.
We claim that (with the notation of Theorem~2.8)
$$(cn^{1/8}\le )\quad \delta_n\le \lambda(n)\le
\sqrt{n}.\leqno (4.13)$$ To verify this, first observe that
$$\lambda(n) =
 \sup\left\{{\|u\|_M\over \|u\|_{\rm min}}\right\}\leqno (4.12)'$$
  {where the supremum} runs over all operator spaces $E,F$
in $OS_n$ and all $u$ in $E\otimes F$. Equivalently, we
have $$\lambda(n) = \sup\{\|i_E\|_M\}\leqno (4.12)''$$
where the supremum runs over all $E$ in $OS_n$ and where $i_E\in
E\otimes
E^*$ represents (as above) the identity on $E$.

\n Indeed, any $u$ as in (4.12)$'$ can be rewritten as $u =
(I_E\otimes \widetilde u)(i_E)$ where $\widetilde u\colon
\ E^*\to F$ is the linear map corresponding to $u$, hence
by (4.1) we have $\|u\|_M \le \|\widetilde u\|_{cb}
\|i_E\|_M = \|u\|_{\rm min} \|i_E\|_M$, and (4.12)$''$
follows.

\n By Corollary~4.3, this implies
$$\lambda(n) = \sup\{d_{SC^*(F_\infty)}(E)\mid E\in OS_n\}.\leqno
(4.12)'''$$
Now since $C^*(F_\infty)$ is separable, with the notation
 of Theorem~2.8,
we clearly have (recall Proposition 2.6 (i)) $\delta_n\le
\sup\{d_{SC^*(F_\infty)}(E)\mid E\in OS_n\}$, hence
the left side of (4.13)
follows from (4.12)$'''$. For the other side,
note that by (4.4) and (4.5) we have for any $E$ in $OS_n$,\
$d_{SC^*(F_\infty)}(E)\le d_{SK}(E)$
and by [P3,Theorem 9.6] this is $\le \sqrt{n}$. \qed

\proclaim Theorem 4.4. Let $A$ be a $C^*$-algebra and let $H=\ell_2$.
The
following are equivalent.
\item{\rm (i)} $$A\otimes_{\min} B(H)=A\otimes_{M}
B(H).$$
\item{\rm (ii)} For any $\eps>0$ and any finite dimensional subspace
$E\subset A$, there is a subspace $\hat E\subset \C$ such that
$d_{cb}(E,\hat E)<1+\eps.$
\item{\rm (iii)} Same as (ii) with $\eps=0$. Equivalently, every finite
dimensional subspace
of $A$ is completely isometric to a subspace of $\C$.

\pf Assume (i). Then for any finite dimensional subspace
$E\subset A$, consider a completely isometric embedding
$E^*\subset B(H)$ and view $E\otimes E^*$ as a subspace of
$A\otimes_{\min} B(H)$. If we apply Proposition 4.2 when $u\in E\otimes
E^*$ represents the
identity
on $E$, we have $\|u\|_{cb}=1$, hence $\|u\|_{M}=1$ and we immediately
obtain (ii).
Conversely, (ii) clearly implies (i) by Proposition 4.2 and (4.8).
 The fact that (i) implies
(iii) follows from  [EH, Theorem 3.2] (this was kindly pointed out to
us by Kirchberg).
Indeed, it suffices to prove (iii) for finite dimensional operator
systems $E\subset A$.
Then, assuming (i) the second condition in
  [EH, Theorem 3.2] must hold by Lemma 0.7 and the short exactness
  property of
the $\max$-tensor product. Therefore, if we represent
$A$ as a quotient of $C^*(F_I)$ for some free group $F_I$, by [EH,
Theorem 3.2] any unital
completely positive map $v\colon E\to A$ has a unital completely
positive lifting
$\hat v \colon E\to C^*(F_I)$. In particular $E$ embeds completely
isometrically into
$C^*(F_I)$, hence into $\C$ (see the proof of Corollary 3.1). This
shows
(i)$\implies $(iii). Finally,  (iii) $\implies$ (ii)  is trivial. \qed

\n {\bf Remark.} The notion appearing in (ii) above is analogous to
that
of ```finite representability" in Banach space theory.

\vfill\eject

\proclaim Theorem 4.5. Let $X$ be an operator space
and let $c\ge1$ be a constant. The following are
equivalent.
\item{\rm (i)} For any operator space $F$,
we have $X\otimes_{\min}F=X\otimes_M F$ and
 $\N{u}_M\leq c\N{u}_{\min}$ for any $u$ in $X\otimes F$.
\item{\rm (ii)} The same as {\rm (i)} with $F=B(\ell_2)$.
\item{\rm (iii)} For any finite dimensional subspace $E\subset X$ we
have
$$d_{SC^*(F_\infty)}(E)\leq c.$$

\pf (i)$\Rightarrow$(ii) is trivial and (ii)$\Rightarrow$(iii) is clear
by
Proposition 4.2, taking $F=E^*$ and $u\in X\otimes E^*$ associated to
the
inclusion $E\subset X$. Finally assume (iii). Consider
$u$ in
$X\otimes F$. We have
$u\in E\otimes F$ for some finite dimensional subspace $E\subset X$. By
(iii),
for each $\eps>0$, there is a subspace $\tilde E\subset C^*(F_\infty)$
such that
$d_{cb}(E,\tilde E)\leq c+\eps$. By Proposition 4.2,
this implies
$\|u\|_M\le (c+\eps)\|u\|_{\min}$, whence (i). \qed

\n{\bf Remark.} It can be shown that the operator Hilbert space
$OH$ introduced in  [P3] satisfies the equivalent conditions
in  Theorem 4.5 for  $c=1$. In particular we have
$$d_{SC^*(F_\infty)}(OH)=1.$$
This is a consequence of some unpublished
 work by U.~Haagerup, namely
 the inequality (4.14) below
(itself a consequence of [P3, Corollary 2.7]). To explain this
inequality, let $\overline{B(H)}$ be the complex conjugate
of $B(H)$, {\it i.e.} the same space but with the
conjugate complex multiplication. We denote by $x\to
\overline x$ the canonical
anti-isomorphism between $B(H)$ and $\overline{B(H)}$.
Note that $\overline{B(H)}\approx B(\overline H) \approx
B(H)$. In the sequel, we simply denote by $\|\ \|_{\rm max}$
(resp.\ $\|\
\|_{\rm min}$) the max-norm (resp.\ the min-norm) on the
space $B(H)\otimes \overline{B(H)}$.

Then, for any $x_1,\ldots, x_n$ and $y_1,\ldots, y_n$ in $B(H)$ we have
$$\big\|\sum x_i\otimes \overline y_i\big\|_{\rm max} \le \big\| \sum
x_i\otimes \overline x_i\big\|^{1/2}_{\rm min} \big\|\sum y_i\otimes
\overline y_i\big\|^{1/2}_{\rm min}.\leqno (4.14)$$
To check (4.14), we first recall an entirely elementary fact:\ for any
$a_1,\ldots, a_n$, $b_1,\ldots, b_n$ in a $C^*$-algebra $A$, we have
$$\big\|\sum a_ib_i\big\| \le \big\|\sum a_ia^*_i\big\|^{1/2}
\big\|\sum
b^*_ib_i\big\|^{1/2}.$$
Hence if $a_ib_i = b_ia_i$ we also have
$$\big\|\sum a_ib_i\big\| \le \big\|\sum a^*_ia_i\big\|^{1/2}
\big\|\sum
b_ib^*_i\big\|^{1/2}.$$
Applying these inequalities in the case $a_i = x_i\otimes 1$, $b_i =
1\otimes \overline y_i$ we get
$$\leqalignno{\big\|\sum x_i\otimes \overline y_i\big\|_{\rm max} &\le
\big\|\sum x_ix^*_i\big\|^{1/2} \big\|\sum
y^*_iy_i\big\|^{1/2}&(4.14.0)\cr
\big\|\sum x_i\otimes \overline y_i\big\|_{\rm max} &\le \big\|\sum
x^*_ix_i\big\|^{1/2} \big\|\sum y_iy^*_i\big\|^{1/2}.&(4.14.1)}$$

Let us denote  by $A_0$ (resp.\ $A_1$) the space $B(H)^n$ equipped with
the
norm
$\|(x_i)\| = \big\|\sum x^*_ix_i\big\|^{1/2}$ (resp.\ $\|\sum
x_ix^*_i\|^{1/2}$). Moreover for any $(x_i)$ in $B(H)^n$, we denote by
$\|(x_i)\|_{1\over 2}$ the norm in the complex interpolation space
$(A_0,A_1)_{1\over 2}$. Then, by the complex interpolation
theorem applied to the sesquilinear map
$$(x_i), (y_i) \to \sum x_i\otimes \overline y_i,$$
(4.14.0) and (4.14.1) imply that we have
$$\big\|\sum x_i\otimes \overline y_i\big\|_{\rm max} \le
\|(x_i)\|_{1\over
2} \|(y_i)\|_{1\over 2}.\leqno (4.15)$$
On the other hand, by Corollary~2.7 in [P3] we have for all $(x_i)$ in
$B(H)^n$
$$\|(x_i)\|_{1\over 2} = \big\|\sum x_i\otimes \overline
x_i\big\|^{1/2}_{\rm min}.\leqno (4.16)$$
Hence (4.15) implies (4.14).

Finally, let $(T_i)_{i\ge 1}$ be an orthonormal basis in the operator
Hilbert space $OH$ introduced in [P3]. We may assume $OH\subset B(H)$
with
(say) $H=\ell_2$. Recall (see [P3]) that $\left\|\sum\limits^n_1
T_i\otimes
\overline T_i\right\|_{\rm min} = 1$ and for any $x_1,\ldots, x_n$ in
$B(H)$ we have
$$\big\|\sum T_i\otimes x_i\big\|_{\rm min} = \big\|\sum x_i\otimes
\overline x_i\big\|^{1/2}_{\rm min}.$$
Therefore, by (4.14) for any $x_1,\ldots, x_n$ in $B(H)$ we have
$$\big\|\sum T_i\otimes x_i\big\|_{\rm max} \le \big\|\sum x_i\otimes
\overline x_i\big\|^{1/2}_{\rm min} = \big\|\sum T_i\otimes
x_i\big\|_{\rm
min}.$$
Equivalently, we  conclude that $\|\ \|_{\rm max}$ and $\|\ \|_{\rm
min}$
coincide on $OH\otimes B(H)$ so that $X=OH$ satisfies the equivalent
properties in Theorem~4.5.
Note that if $x_i=y_i$ in (4.14) we have
$$\big\|\sum x_i\otimes \overline x_i\big\|_{\rm max} \le \big\|\sum
x_i\otimes \overline x_i\big\|_{\rm min}.$$
In other words, we obtain that
$\|\ \|_{\rm max}$ and $\|\ \|_{\rm min}$ coincide on the ``positive''
cone
in $B(H)\otimes \overline{B(H)}$ formed of all the tensors of the form
$\sum\limits^n_1 x_i\otimes \overline x_i$.

In [P5], a modified version of the identity (4.16) is
proved with an arbitrary semi-finite von~Neumann algebra
in the place of $B(H)$. We refer the reader to a possibly
forthcoming   paper by U.~Haagerup and the second author
for extended
versions of (4.14) and (4.15).\bigskip

\vfill\eject

 \centerline{\bf References}

\item{[AO]}   C. Akeman and P. Ostrand.
 Computing norms in group $C^*$-algebras.
Amer. J. Math. 98 (1976), 1015-1047.

\item{[Ba]} R. Baire. Sur les fonctions des variables r\'eelles.
 Ann. di Mat. (3) 3 (1899) 1-123.

 \item{[B1]} D. Blecher. Tensor products of
 operator spaces II. (Preprint)
 1990.  Canadian J. Math. 44 (1992) 75-90.

 \item{[B2]}  $\underline{\hskip1.5in}$. The standard dual of an
 operator space.
  Pacific J. Math. 153 (1992) 15-30.

 \item{[B3]}  $\underline{\hskip1.5in}$.
 Tracially completely bounded
multilinear maps on $C^*$-algebras.
Journal of the London Mathematical Society,  39 (1989) 514-524.

 \item{[BP]} D. Blecher and V. Paulsen. Tensor products of operator
 spaces.
 J. Funct. Anal. 99 (1991) 262-292.

\item{ [DCH]} J. de Canni\`ere   and  U. Haagerup.
Multipliers of the Fourier algebras of some simple Lie
groups and their discrete subgroups.  Amer. J. Math.   107
(1985), 455-500.

\item{ [DHV]} P. de la Harpe and A. Valette. La propri\'et\'e
T de Kazhdan pour les groupes localement compacts.
Ast\'erisque, Soc. Math. France 175 (1989).

 \item{[EH]} E. Effros and U. Haagerup. Lifting problems
and local reflexivity for $C^*$-algebras. Duke Math. J.
52 (1985) 103-128.

 \item{[ER1]} $\underline{\hskip1.5in}$. A new approach
to operator spaces.
 Canadian Math. Bull.
34 (1991) 329-337.

 \item{[ER2]} $\underline{\hskip1.5in}$. On the abstract
characterization of
 operator spaces. Proc. Amer. Math. Soc. 119 (1993) 579-584.

\item{[H1]} U. Haagerup. The Grothendieck inequality for bilinear forms
on
$C^*$-algebras. Advances in Math. 56 (1985) 93-116.

\item{[H2]} $\underline{\hskip1.5in}$. An example of a non-nuclear
 $C^*$-algebra which
 has the metric approximation property. Invent. Math. 50 (1979)
 279-293.

\item{[H3]} $\underline{\hskip1.5in}$. Injectivity and decomposition of
completely
 bounded maps in ``Operator algebras and their connection with Topology
 and
 Ergodic Theory''. Springer Lecture Notes in Math. 1132
(1985) 170-222.

 \item{[HP]} U. Haagerup and G. Pisier.  Bounded linear
operators between
 $C^*$-algebras. Duke Math. J. 71 (1993) 889-925.

 \item{[I]} T. Itoh. On the completely bounded maps
 of a $C^*$-algebra to its dual space.
Bull. London math. Soc. 19 (1987) 546-550.

 \item{[Ki1]} E. Kirchberg. On subalgebras of the
CAR-algebra. (Preprint, Heidelberg, 1990)
J. Funct. Anal. (To appear)

 \item{[Ki2]} $\underline{\hskip1.5in}$. On non-semisplit extensions,
 tensor products and
exactness of group $C^*$-algebras.
Invent. Math. 112 (1993) 449-489.

 \item{[Ki3]}$\underline{\hskip1.5in}$. Commutants of unitaries in UHF
 algebras
and functorial properties of exactness. J.
 reine angew. Math. To appear.

 \item{[Kr]} J. Kraus. The slice map problem and
approximation properties. J. Funct. Anal. 102 (1991)
116-155.

\item{[Ku]} W. Kuratowski. Topology. Vol. 1. (New edition
translated from the french.) Academic Press, New-York 1966.

 \item{[Kw]} S. Kwapie\'n.  On operators factorizable
through $L_p$-spaces.
 Bull. Soc. Math. France, M\'emoire
31-32,(1972) 215-225.

 \item{[La]} C. Lance. On nuclear $C^*$-algebras. J. Funct. Anal.
12 (1973) 157-176.

 \item{[Pa1]} V. Paulsen.  Completely
bounded maps and dilations. Pitman Research Notes 146.
Pitman Longman (Wiley) 1986.

\item{[Pa2]} $\underline{\hskip1.5in}$. Representation of
Function algebras, Abstract
 operator spaces and Banach space Geometry. J.
Funct. Anal. 109 (1992) 113-129.

\item{[Pa3]} $\underline{\hskip1.5in}$. The maximal operator
space of a normed space. To appear.

 \item{[P1]} G. Pisier. Exact operator
spaces. Colloque sur les alg\`ebres d'op\'erateurs (Orl\'eans 1992).
 Ast\'erisque. Soc. Math. France. To appear.

\item {[P2]} $\underline{\hskip1.5in}$. Factorization of linear
operators and the Geometry of Banach spaces.  CBMS
(Regional conferences of the A.M.S.)    60, (1986),
Reprinted with corrections 1987.

\item{[P3]} $\underline{\hskip1.5in}$. The operator Hilbert space $OH$,
complex interpolation and tensor norms.
Memoirs Amer.
Math. Soc. (submitted)

 \item{[P4]} $\underline{\hskip1.5in}$. Factorization of operator
 valued analytic functions.
 Advances in Math. 93 (1992) 61-125.

\item{[P5]} $\underline{\hskip1.5in}$. Projections from
a von Neumann algebra
onto a subalgebra. Bull. Soc. Math. France. To appear.

\item{[Ru]} Z. J. Ruan. Subspaces of $C^*$-algebras. J. Funct. Anal. 76
 (1988) 217-230.

\item{[Sa]} S. Sakai. $C^*$-algebras and $W^*$-algebras.
Springer Verlag New-York, 1971.

\item{[Sm]} R. R. Smith. Completely bounded maps
between $C^*$-algebras. J. London Math. Soc. 27
(1983) 157-166.

\item{[Ta]} M. Takesaki. Theory of Operator Algebras I.
Springer-Verlag New-York 1979.

\item{[Tr]} S. Trott. A pair of generators for the unimodular group.
Canad. Math. Bull. 3 (1962) 245-252.

\item{[V]} D. Voiculescu. Property T
and approximation of operators.
Bull. London Math. Soc. 22 (1990) 25-30.

\item{[VDN]} D. Voiculescu, K. Dykema, A. Nica. Free
random variables. CRM Monograph Series,
Vol. 1,  Amer. Math. Soc., Providence RI.

\item{[W1]} S. Wassermann. On tensor products of certain
group $C^*$-algebras. J. Funct. Anal. 23 (1976) 239-254.

\item{[W2]} $\underline{\hskip1.5in}$.
 Exact $C^*$-algebras and related topics.
Lecture notes series 19, Seoul National University, 1994.

\vskip12pt

\vskip12pt

M. Junge:

Mathematisches Seminar

CAU  Kiel,

24098 Kiel, Germany

\vskip12pt
G. Pisier:

Texas A\&M  University

College Station, TX 77843, U. S. A.

and

Universit\'e Paris 6

Equipe d'Analyse, Bo\^\i te 186,

75252 Paris Cedex 05, France
\end

Assume $F\subset B({\cal H})$. By Kirchberg's
(4.3) we have
$\N{~}_{\min}=\N{~}_{\max}$ on
$C^*(F_\infty)\otimes B({\cal H})$, hence a fortiori by (4.8),
$\N{~}_{\min}=\N{~}_M$ on $\tilde E\otimes
F$ which immediately yields (i) using (4.1) and $d_{cb}(E,\tilde E)\leq
c+\eps$.

\n {\bf Remark.} Let $\vp_n$ be the infimum
of the set of numbers $\vp$ such that there is a countable
$\vp$-net in $OS_n$. By [P1]
we have $d_{SK}(E)\le n^{1/2}$ for all $E$ in $OS_n$.
Therefore $\vp_n\le n^{1/2}$. It would be interesting to
estimate the asymptotic behaviour of $\vp_n$ when $n\to \infty$.